\documentclass [smallcondensed] {svjour3}

\usepackage{xcolor}
\usepackage{graphicx}
\usepackage{caption}

\usepackage{changes,cancel}

\newcommand{\pxx}{p_{xx}}

\newcommand{\pyy}{p_{yy}}
\newcommand{\pxy}{p_{xy}}
\newcommand{\N}{\mathcal{N}}
\newcommand{\IS}{\beta}

\usepackage{amsmath,amsfonts}
\usepackage{hyperref}
\newcommand{\R}{\mathbb{R}}

\newcommand{\D}{\mathrm{d}}
\newcommand{\DT}{{\Delta}t}

\newcommand{\pder}[2]{\partial_{#2}{#1}}

\newcommand{\tder}[1]{\frac{\mathrm{d}}{\mathrm{d}#1}}

\newcommand{\diverg}{\nabla\cdot}

\newcommand{\Popt}{P_{\text{OPT}}}

\newcommand{\Sref}{S_{\text{ref}}}
\newcommand{\Scoa}{S_{\text{coa}}}

\begin{document}
\title{Adaptive Mesh Refinement for Hyperbolic Systems based on Third-Order Compact WENO Reconstruction}
\titlerunning{AMR based on Third-Order Compact WENO Reconstruction}
\author{M. Semplice \and A. Coco \and G. Russo }
\institute{
M. Semplice
\at
Dipartimento di Matematica ``G. Peano'' -
Universit\`a di Torino -
Via C. Alberto, 10 -
10123 Torino (Italy)
\email{matteo.semplice@unito.it}
\and
A. Coco
\at
Bristol University, 
Queen’s Road, 
Bristol BS8 1RJ 
(United Kingdom)
\email{Armando.Coco@bristol.ac.uk}
\and 
G. Russo
\at
Dipartimento di Matematica -
Universit\`a di Catania -
Viale Andrea Doria, 6 -
95125 Catania (Italy)
\email{grusso@unict.it}
}
\date{Received: date / Accepted: date (The correct dates will be entered by the editor)}

\maketitle

\begin{abstract}
In this paper we generalize to non-uniform grids of quad-tree type the Compact WENO reconstruction of Levy, Puppo and Russo (SIAM J. Sci. Comput., 2001), thus obtaining a truly two-dimensional non-oscillatory third order reconstruction with a very compact stencil and that does not involve mesh-dependent coefficients. This latter characteristic is quite valuable for its use in h-adaptive numerical schemes, since in such schemes the coefficients that depend on the disposition and sizes of the neighboring cells (and that are present in many existing WENO-like reconstructions) would need to be recomputed after every mesh adaption.

In the second part of the paper we propose a third order h-adaptive scheme with the above-mentioned reconstruction, an explicit third order TVD Runge-Kutta scheme and the entropy production error indicator proposed by Puppo and Semplice (Commun. Comput. Phys., 2011). After devising some heuristics on the choice of the parameters controlling the mesh adaption, we demonstrate with many numerical tests that the scheme can compute numerical solution whose error decays as $\langle N\rangle^{-3}$, where $\langle N\rangle$ is the average number of cells used during the computation, even in the presence of shock waves, by making a very effective use of h-adaptivity and the proposed third order reconstruction.

\end{abstract}

\section{Introduction}

We consider numerical approximations of the solution of an initial value problem for hyperbolic
system of $m$ conservation laws in a domain $\Omega\subset\R^d$
\begin{equation}
\label{eq:pde}
\pder{u}{t} + \diverg f(u) = 0
\end{equation}
with suitable boundary conditions. The finite volume formulation of \eqref{eq:pde} reads
\begin{equation}
\label{eq:finvol}
\tder{t} \int_{\Omega_j} u(t,x)\D{x}
+ \oint_{\partial \Omega_j} f(u(t,x))\cdot\vec{n}\:\D{l}
=0,
\end{equation}
where $\left\{ \Omega_1, \Omega_2, \ldots, \Omega_N \right\}$ denotes a partition of $\Omega$.
A semidiscrete finite volume scheme tracks the evolution in time of the cell
averages
\[ U_j = \frac{1}{|\Omega_j|}\int_{\Omega_j} u(t,x)\:\D{x} \] and, in
order to approximate the boundary integrals in \eqref{eq:finvol} one
would need to know the point values of the approximate solution at
suitable quadrature points along the boundary of each cell
$\Omega_j$. This is accomplished by a reconstruction procedure that
computes point values of $u$ at quadrature nodes along $\partial
\Omega_j$ from the cell average $U_j$ and the cell averages of the
neighbors of $\Omega_j$.

In order to obtain a convergent scheme, the reconstruction procedure
must satisfy requirements of mass conservation, accuracy (on smooth
data) and have some non-oscillatory property, controlling the
boundedness of the total variation (on non-smooth data).

A robust method that yields high order accuracy and possesses
non-oscillatory properties is the so-called WENO reconstruction. In
its original formulation \cite{Shu97}, in the case of a uniform grid
in one space dimension, one starts by considering all the possible $k+1$
stencils (from downwind, to central, to upwind) of width $k+1$
containing the cell $\Omega_j$ and the corresponding polynomials
$p_1(x),\ldots,p_k(x)$ of degree $k$ that interpolate the data in the sense of
the cell
averages on each of the stencils. One can also determine the so-called
linear weights $\alpha_\gamma^\pm$ such that $\sum_\gamma\alpha_\gamma^\pm
p_\gamma(x)$ yields $2k+1$ order accuracy when evaluated at cell
boundaries $x_{j\pm1/2}$ (that is the accuracy given by interpolating in the sense of the cell
averages $U_{j-k},\ldots,U_{j+k}$ with a single polynomial of degree
$2k$). The reconstruction is then taken to be a convex linear
combination $p(x)=\sum_\gamma \tilde{\alpha}_\gamma p_\gamma(x)$,
where the nonlinear weights $\tilde{\alpha}_\gamma$ are chosen
according to the smoothness of the data in the $\gamma$-th stencil
and in such a way that, $\tilde{\alpha}_\gamma \simeq\alpha_\gamma$
for all values of $\gamma$ if the solution is smooth in the union of
all the stencils and $\tilde{\alpha}_\gamma \simeq0$ if the solution
is not smooth in the $\gamma$-th stencil.  Note that there will be
two different sets of linear weights $\left\{\alpha_\gamma^+\right\}$ and $\left\{\alpha_\gamma^-\right\}$ respectively for the right- and left- boundary of
the cell $\Omega_j=[x_{j-1/2},x_{j+1/2}]$, but this can be handled
efficiently in one space dimensions or even in higher dimensional
setting if one has a Cartesian grid and employs dimensional splitting~\cite{Shu97}.

The situation is far more complicated for non-Cartesian grids,
especially for higher order schemes. In fact, in two or more
dimensions, in order to approximate the boundary integral in
\eqref{eq:finvol}, one uses a suitable quadrature formula for each
face of $\partial\Omega_j$. Up to schemes of order $2$, the midpoint
rule can be employed, but for schemes of higher order, one has to use
a quadrature formula with more than one node per face and this
increases the cardinality of the sets of linear weights that should be used.

One approach, typical of the ClawPack software, is to build the computational grids from
regular cartesian patches. We point out that ClawPack now includes higher order solvers (SharpClaw) \cite{Ketcheson2011}, but their code has been so far released only in one space dimension \cite{pyclaw}. The more traditional approach is to employ a single grid composed of cells of different sizes. In this framework, WENO schemes on arbitrary triangular meshes were
constructed in \cite{HuShu:1999}. For efficiency, the third order
reconstruction suitable for a two-point quadrature formula on each
edge, requires the precomputation and storage of $90$ linear weights
and $180$ other constants for the smoothness indicators for each
triangle in the mesh. Moreover, upon refinement or coarsening, one
would have to recompute all the above mentioned constants in all cells
that include the newly-created one in their stencil (usually 10
cells). Additionally, in the general case, the linear weights are not
guaranteed to be positive: see \cite{ShiHuShu:2002} for the adverse
effects of this non-positivity and for ways to circumvent this.
Similar issues are expected to arise when using locally refined
quadrangular grids, like those of Fig.\ \ref{fig:grid2d} as there are many
possible configurations for the number, size and relative position of the
neighbors. More recently \cite{DK:2007:linear} considers a WENO
scheme that combines central and directional polynomials of degree
$k$, requiring very large stencils.

Many of these difficulties are linked to the otherwise desirable aim
of achieving enhanced accuracy of order $2k+1$ at the quadrature
nodes. Relaxing this requirement, a reconstruction of order $k+1$ can
obviously be achieved by combining polynomials of order up to $k$ on
each stencil using a single set of linear weights for every quadrature
node of every face and for all disposition of cells in the
neighborhood. In fact, this procedure gives a $k$-degree polynomial
that represents a uniform approximation of order $k+1$ on the cell and
thus can just be evaluated at all quadrature nodes. The accuracy would
thus be similar to the accuracy of a ENO procedure, but using a
nontrivial linear combination of all polynomials avoids the
difficulties of ENO reconstructions on very flat solutions~\cite{RoMe:90}.

Another problem arising is related to the stencil widths. In fact, in
order to compute a polynomial of degree $k$ interpolating the cell
averages of the numerical solution using a fully upwind stencil in
some direction, one needs to use data that is quite far away from the
cell $\Omega_j$. This poses a question on how efficiently one can
gather the information of the local topology of the triangulation and,
since this may change at every timestep in an adaptive mesh refinement
code, it would be a clear advantage to employ a reconstruction that
can achieve order $k+1$ on smooth solutions by using a compact stencil.

A step in this direction was set by \cite{LPR:2001}, that introduced the Compact WENO (CWENO)
reconstruction of order $3$ for regular Cartesian grids. The
main idea is to consider a central stencil composed only by the cells
that share at least a vertex with $\Omega_j$ and to construct an ``optimal''
reconstruction polynomial $\Popt$ of degree $2$ that interpolates
exactly the cell averages in this stencil. In order to control the
total variation increase possibly caused by the presence of a
discontinuity in the stencil, this polynomial is decomposed as
\[
\Popt = \alpha_0 P_0 + \sum_{\gamma=1}^{D(d)} \alpha_\gamma P_\gamma ,
\] 
where the $P_\gamma$'s for $\gamma>0$ constitute a finite set of $D(d)=2^d$
first degree polynomials interpolating in the sense of the cell averages in a
directional sector of the central stencil in the direction of each of the $2^d$ vertices of the $d-$dimensional cubic cell.
$P_0$ is a second
degree polynomial computed by difference from $\Popt$. The
reconstruction polynomial is then given by
\[ 
P = \tilde{\alpha}_0 P_0 + \sum_{\gamma=1}^{D(d)} \tilde{\alpha}_\gamma P_\gamma,
\] 
where the coefficients $\tilde{\alpha}_{\gamma}$ are computed from the
linear weights $\alpha_{\gamma}$ as usual, with the help of smoothness
indicators. This procedure allows to obtain third order accuracy in
smooth regions, degrades to a linear reconstruction near
discontinuities, but employs a very compact stencil.
Notice that a natural choice is to take the $\alpha_\gamma = (1-\alpha_0)/D(d)$ for each $\gamma = 1, \dots, D(d)$. Therefore the set of linear weights depends on a single parameter $\alpha_0 \in (0,1)$. In practice we use $\alpha_0=1/2$ in the whole paper. In our experience, the use of higher values of $\alpha_0$ are beneficial only in problems without discontinuities. 

This reconstruction does not aim at computing directly the point values at specific
quadrature points on the cell boundary, but gives a single
reconstruction polynomial that can then be evaluated at any point in the
cell with uniform accuracy. The degree of accuracy is of course the same than the
one achievable by the standard WENO technique with the same stencils,
but this approach seems more suitable in the context of
non-uniform and adaptive grids, because it requires a single set of
linear weights. This feature turns out to be quite useful in the case
of balance laws with source terms, where the reconstruction at
quadrature points inside the cell is required for the well-balanced
quadrature of the source term, see \cite{PS:shentropy}.

In order to construct a third order scheme suitable for an adaptive
mesh refinement setting, in this paper we extend the idea of the CWENO
reconstruction to nonuniform meshes. In one space dimension the
extension is straightforward, while in higher dimensions complications
arise due to the variable cardinality of the set of first neighbors
of a given cell. These are solved in Section 2 by considering an
interpolation procedure that enforces exact interpolation of the cell
average in the central cell but only a best fit in the least
squares sense of the other cell averages in the stencil. Section 3 is
devoted to analyzing the behavior of an adaptive scheme based on some
error indicator and the employment of cells of sizes between $H$ and
$H/2^L$, in order to provide guidance on the number of ``levels'' $L$
that need to be used and on the response of the scheme to the choice
of the refinement threshold. The case of global timestepping is considered 
in this paper. Finally, Section 4 presents numerical
tests in one and two space dimensions that demonstrate both the
properties of the reconstruction introduced in Section 2 and
corroborate the ideas on the behavior of adaptive schemes illustrated
in Section 3. Section 5 contains the conclusions and perspectives for
future work.

\section{Third order, compact stencil reconstruction}

In our approach, we try to combine the simplicity of the Cartesian meshes, the compactness of the CWENO reconstruction and the flexibility of the quad-tree local grid refinement for the construction of h-adaptive schemes.
We thus restrict ourselves to grids where each cell $\Omega_j$ is a $d-$dimensional cube of edge $h_j$
and may be refined only by splitting it in $2^d$ equal
parts. Such grids are saved in binary/quad-/oct-trees.
Let us introduce some notation. Let $\Omega$ be the entire domain, $\left\{ \Omega_j, j=1,\ldots,N \right\}$ the set of all cells, $\vec{x}_j$ the center of cell $\Omega_j$, and $|\Omega_j| = h_j^{d}$.

\subsection{One space dimension}

In one space dimension the reconstruction is a straightforward
extension of the reconstruction described in~\cite{LPR:2001} to the
case of non-uniform grids.  Let $u(x)$ be a function defined in
$\Omega=[a,b]$ and let us suppose we know the mean values of $u$ in each
cell $\Omega_j$, i.e.\ we know
\[
 U_j=\frac{1}{h_j} \int_{\Omega_j} u(x)dx, \qquad j=1, \ldots,N,
\]
with $h_j=\left|\Omega_j\right|$.
Now let us fix a cell $\Omega_j$. In order to compute the
reconstruction in that cell, we use a compact stencil, namely we use
only the mean values of $u$ in that cell, $U_j$, and in the first
neighbors, $U_{j-1}$ and $U_{j+1}$. If $u(x)$ is locally smooth,
%in $\displaystyle \cup_{i=-1}^1 \Omega_{j+i}$ 
one may choose the third order accurate reconstruction given by the
parabola that interpolates the data $U_{j+i}$ in the sense of
cell-averages to enforce conservation. We shall call such a parabola the {\em optimal}
polynomial $\Popt(x)$:
\begin{equation*}
\frac{1}{h_{j+i}} \int_{\Omega_{j+i}} \Popt(x) dx = U_{j+i}, \qquad i=-1,0,1.
\end{equation*}
The optimal polynomial $\Popt(x)$ is completely determined by these
conditions, and its expression is given by:
\[
 \Popt(x) = U_j + p_x (x-x_j) + \frac{1}{2} p_{xx} \left((x-x_j)^2-\frac{h_j}{12} \right),
\]
with
\[
 p_x=\frac{(h_j+2h_{j-1}) U[j-1;j]+(h_j+2h_{j+1}) U[j;j+1]}{2(h_{j-1}+h_j+h_{j+1})},
\]
\[
 p_{xx}=\frac{3(2h_j+h_{j-1}+h_{j+1})U[j-1;j;j+1]}{2(h_{j-1}+h_j+h_{j+1})},
\]
where
\[
\]
\[
 U[{j-1};{j}]=\frac{U_{j}-U_{j-1}}{x_{j}-x_{j-1}}, \qquad U[{j};{j+1}]=\frac{U_{j+1}-U_j}{x_{j+1}-x_j}
\]
\[
U[{j-1};{j};{j+1}]=\frac{U[{j};{j+1}]-U[{j-1};{j}]}{x_{j+1}-x_{j-1}}.
\]
If $u(x)$ is not smooth in $\displaystyle \cup_{i=-1}^1\Omega_{j+i}$,
then this reconstruction would be oscillatory. Following the idea
of~\cite{LPR:2001}, we compute two linear functions $P_\gamma$,
$\gamma=1,2$ in such a way $P_\gamma$ matches the cell averages $U_j$
and $U_{j+2\gamma-3}$, namely:
\[
 P_\gamma(x)=U_j+U[{j-2+\gamma};{j-1+\gamma}](x-x_j), \qquad \gamma=1,2,
\]
and a parabola $P_0$ determined by the relation:
\[
 \Popt = \alpha_0 P_0 + \sum_{\gamma=1}^2 \alpha_\gamma P_\gamma,
\]
where the coefficients $\alpha_\gamma$ can be chosen arbitrarily,
provided $\alpha_\gamma>0$ and $\sum_{\gamma=0}^2
\alpha_\gamma=1$. 
%The choice is arbitrary since every set of coefficients lead to use $\Popt$ as reconstruction in smooth regions.
In practice we use $\alpha_0=1/2$ and $\alpha_\gamma=1/4$,
$\gamma=1,2$.  The final reconstruction polynomial is:
\begin{equation}\label{finalrec1d}
 P= \tilde{\alpha}_0 P_0 + \sum_{\gamma=1}^2 \tilde{\alpha}_\gamma P_\gamma
\end{equation}
with
\[
 \tilde{\alpha}_\gamma = \frac{\omega_\gamma}{\sum_{\delta=0}^2 \omega_\delta}, \qquad \omega_\gamma = \frac{\alpha_\gamma}{(\epsilon+\IS_\gamma)^2}, \quad \gamma=0,1,2.
\]
The smoothness indicators, $\IS_\gamma$, are responsible for detecting
large gradients or discontinuities and to automatically switch to the
stencil that generates the least oscillatory reconstruction in such
cases. 
On the other hand, if $u(x)$ is smooth the smoothness indicators should be expected to be close to each other and thus $\tilde{\alpha}_\gamma \approx \alpha_\gamma$ and $P(x) \approx \Popt (x)$ for every set of coefficients $\alpha_\gamma$.
Following~\cite{Shu97,LPR:2001}, we define the
smoothness indicators as:
\begin{equation}\label{IS1d}
 \IS_\gamma=\sum_{l=1}^2 \int_{\Omega_j} |\Omega_j|^{2l-1} (P_\gamma^{(l)}(x))^2 dx, \qquad \gamma=0,1,2.
\end{equation}
Let us rewrite the central polynomial $P_0$ as:
\[
 P_0 = U_j + p_x^0 (x-x_j) + \frac{1}{2} p_{xx}^0 \left((x-x_j)^2-\frac{h_j^2}{12} \right).
\]
A direct computation of \eqref{IS1d} yields:
\[
 \IS_0=\frac{13}{12} h_j^4 (p_{xx}^0)^2+h_j^2 (p_x^0)^2,
\]
\[
 \IS_1 = U[{j-1};{j}]^2 h_j^2, \qquad \IS_2 = U[{j};{j+1}]^2 h_j^2.
\]
% The constant $\epsilon$ is chosen in such a way $\alpha_\gamma \approx
% \tilde{\alpha}_\gamma$ in smooth regions (and therefore $P \approx
% \Popt$), while near discontinuities we want that the weights of the
% parabola and of one of the one-sided linear reconstructions will be
% practically zero, and we will be left with the other one-sided linear
% reconstruction. This is accomplished requiring that
% (see~\cite{LPR:2001}):
% \begin{itemize}
%  \item $\epsilon+h_j^2 u_x^2 \gg |IS-h_j^2u_x^2|$ in smooth regions,
% \item $\epsilon\ll IS$ near discontinuities.
% \end{itemize}
% In practice we use $\epsilon=...$. \armando{scrivere dopo i test!}

It is known that the role of $\epsilon$ goes beyond the simple
avoidance of zero denominators in the computation of the nonlinear
weights, but its choice can influence the order of convergence of the
method. In fact, a fixed value of $\epsilon$ is not suitable to
achieve the theoretical order of accuracy both on coarse and fine
grids. Namely, a high value for $\epsilon$ may yield oscillatory
reconstructions, while low $\epsilon$ values yields an estimated order
of convergence that matches the theoretical one only asymptotically
and often only for very fine grid spacings. Among the different
solutions proposed in the literature, the mappings of
\cite{HAP:2005:mappedWENO,FHW:12:mappedWENO} do not apply in a
straightforward way to our compact WENO reconstruction technique,
while taking an $h$-dependent $\epsilon$ as in \cite{ABBM:11:epsilonh}
yields an improvement of the reconstruction, but we found
experimentally that the scaling $\epsilon\propto h$ works best in our
situation (see the numerical tests). Our choice is further supported by the work of Kolb \cite{Kolb14} that analyzed the optimal convergence rate of CWENO schemes depending on the choice of $\epsilon(h)$ on uniform schemes. While on a uniform grid one can choose a constant value of $\epsilon$ that guarantees good numerical results, when using adaprive grids, with cell size $h$ varying of several orders of magnitude, selecting the proper dependence $\epsilon(h)$ is of paramount importance.

\subsection{Higher space dimension}
Although our reconstruction can be performed in any space dimension,
we describe here for simplicity the two-dimensional case, giving hints for its generalization.
We start describing the structure of the adaptive grid.

\subsubsection{Quad-tree grid in 2D}
The adaptive grid is recursively generated starting by a coarse uniform Cartesian mesh of grid size $H$ at level $l=0$.
Then each cell is possibly (according to some criterion) recursively subdivided into four equal squares. At the end of the recursive subdivision, the grid structure is described by a quad-tree.

\begin{figure}[!hbt]
  \begin{minipage}[c]{0.49\textwidth}
    \centering
    \includegraphics[width=0.99\textwidth]{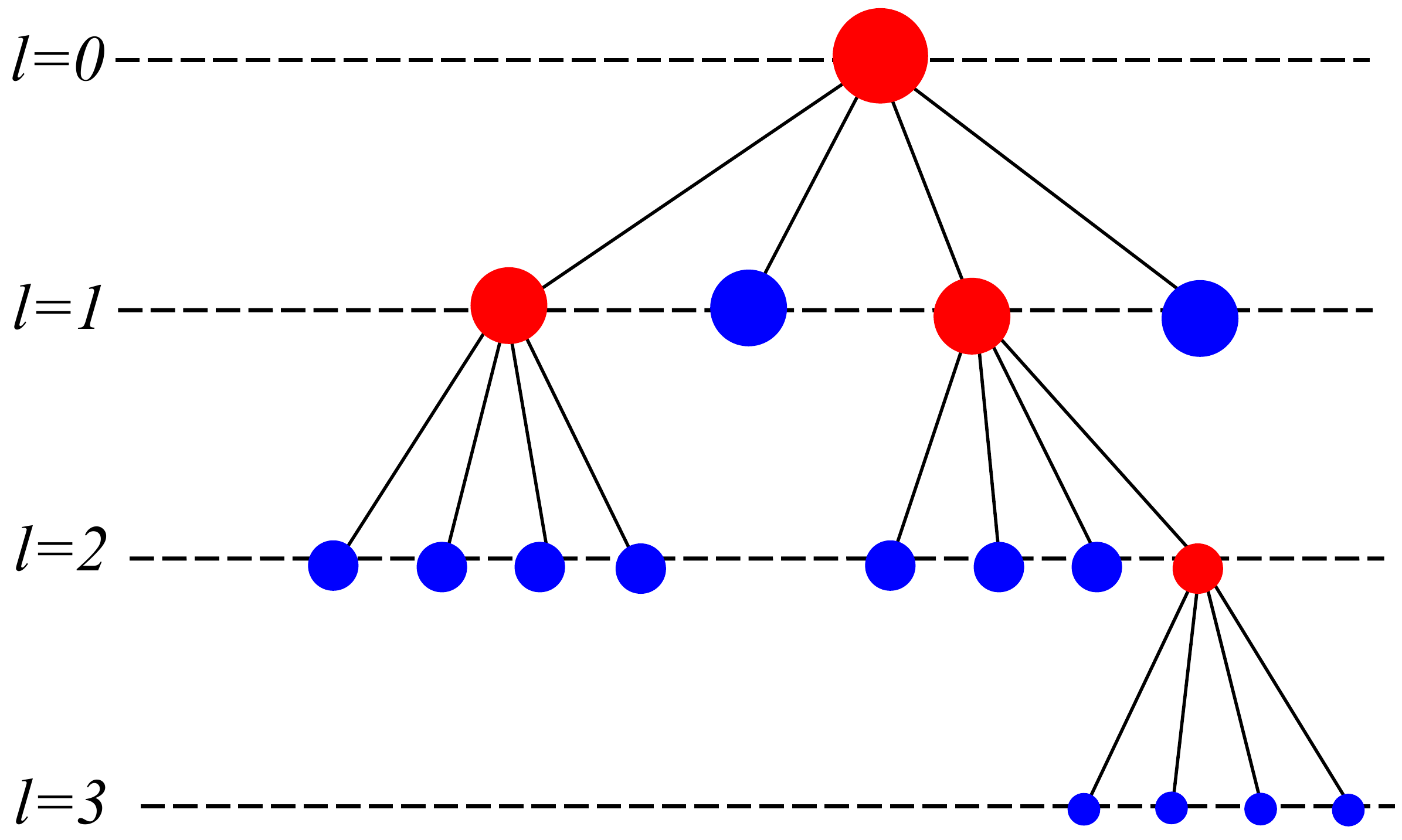}
  \end{minipage}
  \ \hspace{2mm} \hspace{3mm} \
  \begin{minipage}[c]{0.49\textwidth}
    \centering
    \includegraphics[width=0.99\textwidth]{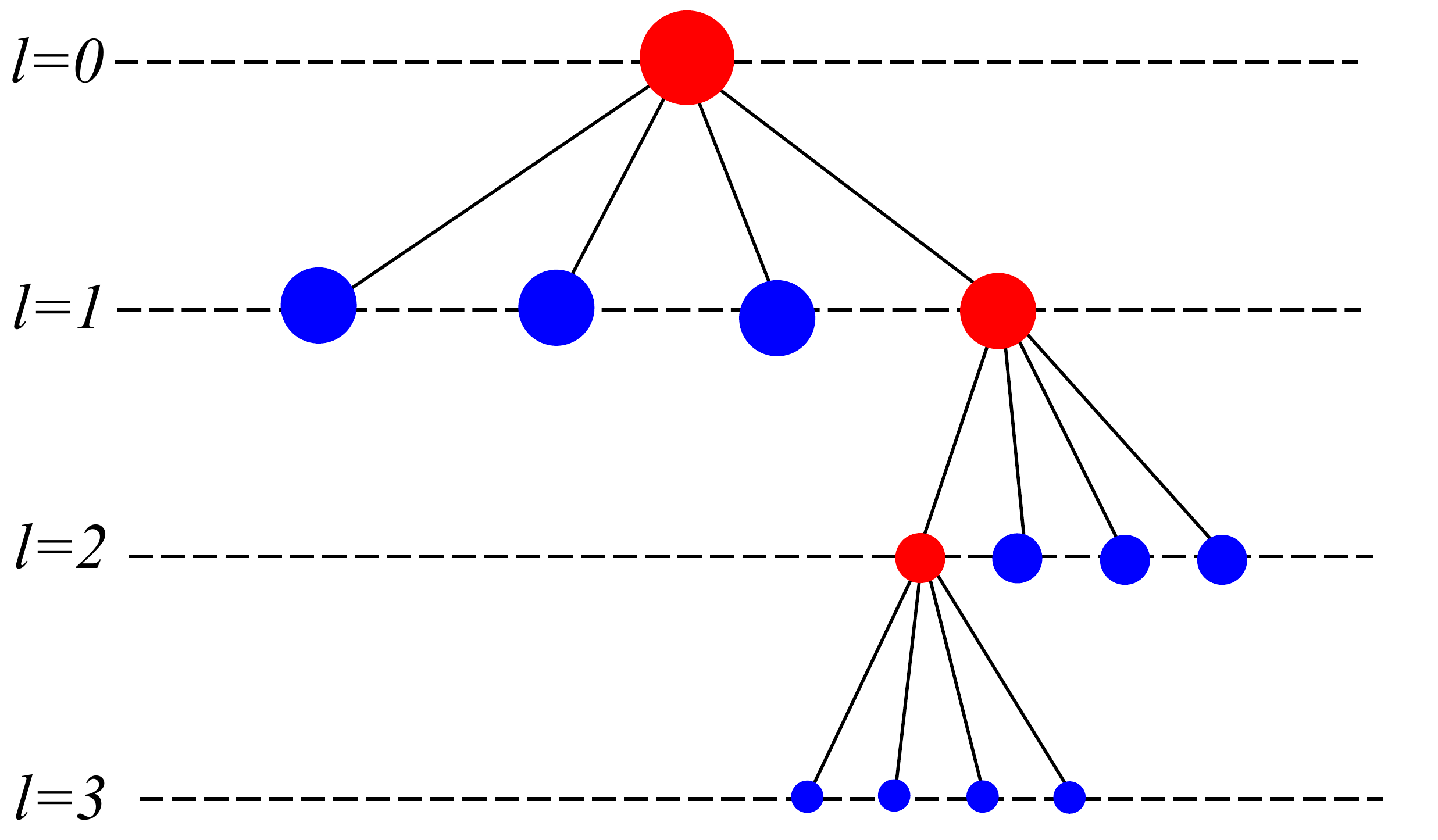}    
  \end{minipage}
    \begin{minipage}[c]{0.49\textwidth}
    \centering
        \includegraphics[width=0.89\textwidth]{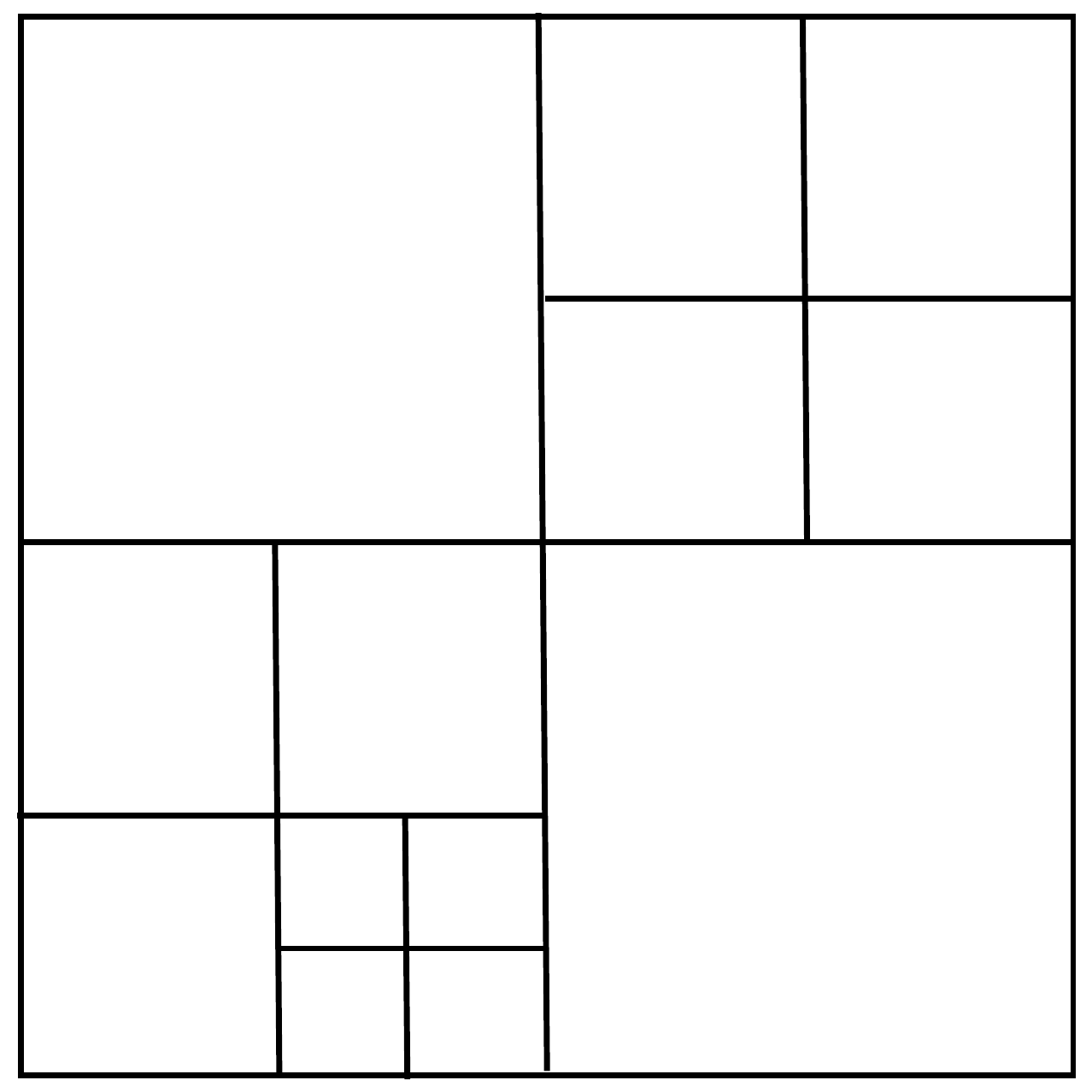}
  \end{minipage}
  \ \hspace{2mm} \hspace{3mm} \
  \begin{minipage}[c]{0.49\textwidth}
    \centering
        \includegraphics[width=0.89\textwidth]{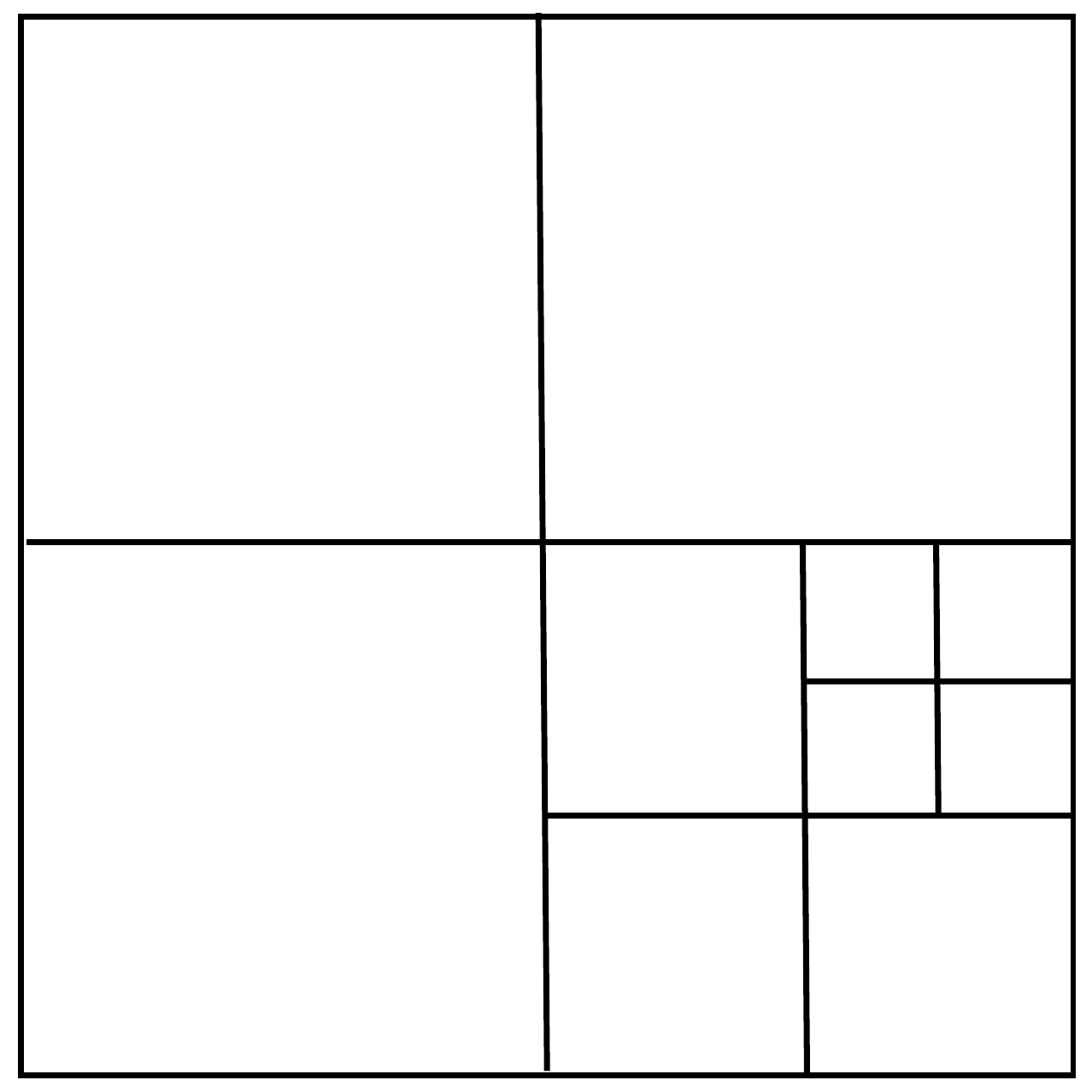}
  \end{minipage}
  \caption{ {Two examples (left and right) of subdivision of an original cell with $L=3$ levels of refinement: quad tree (up) and grid (down). }  }
  \label{fig:quadtree_grid}
\end{figure}

Two examples of subdivision of an original cell with $L=3$ levels of refinement are illustrated in Fig.\ \ref{fig:quadtree_grid}, together with the corresponding quad-tree.
The cell corresponding to the level $l=0$ is the root of the quad-tree. Each cell of level $l \in \left\{ 1, \ldots, L \right\}$ has a father cell, which corresponds to its neighbor node in the quad-tree at level $l-1$. The four nodes connected to the father node are called the children of the node. The children of each subdivided cell are given some prescribed ordering (e.g. counter-clockwise starting from the upper-right as in the Figure).

\subsubsection{CWENO reconstruction in 2D adaptive grid}

\begin{figure}[!hbt]
  \begin{minipage}[c]{0.45\textwidth}
    \centering
    \captionsetup{width=0.75\textwidth}
    \includegraphics[width=0.99\textwidth]{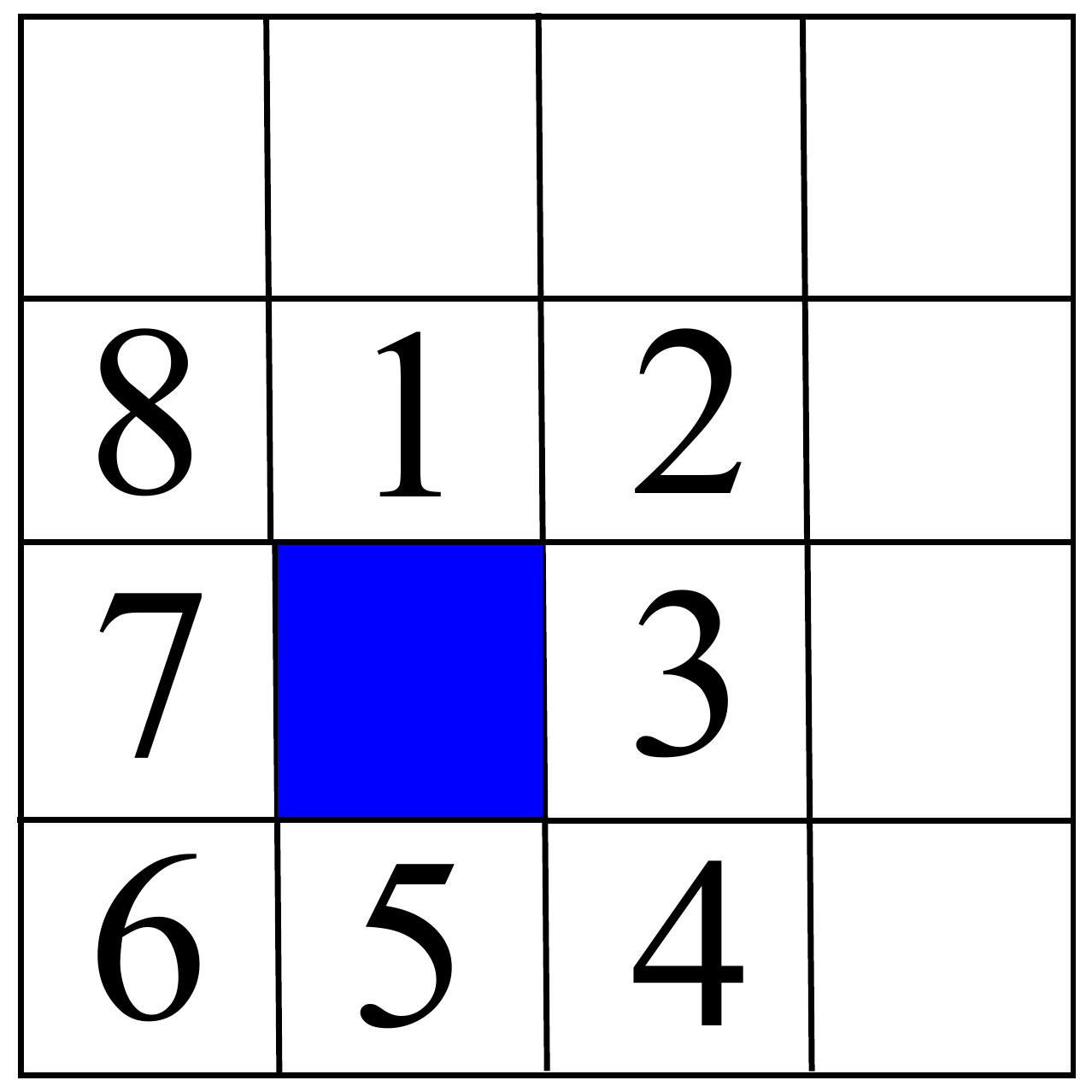}
    % \caption{   }
    % \label{fig:grid2d}
  \end{minipage}
  \ \hspace{2mm} \hspace{3mm} \
  \begin{minipage}[c]{0.45\textwidth}
    \centering
    \captionsetup{width=0.75\textwidth}
    \includegraphics[width=0.99\textwidth]{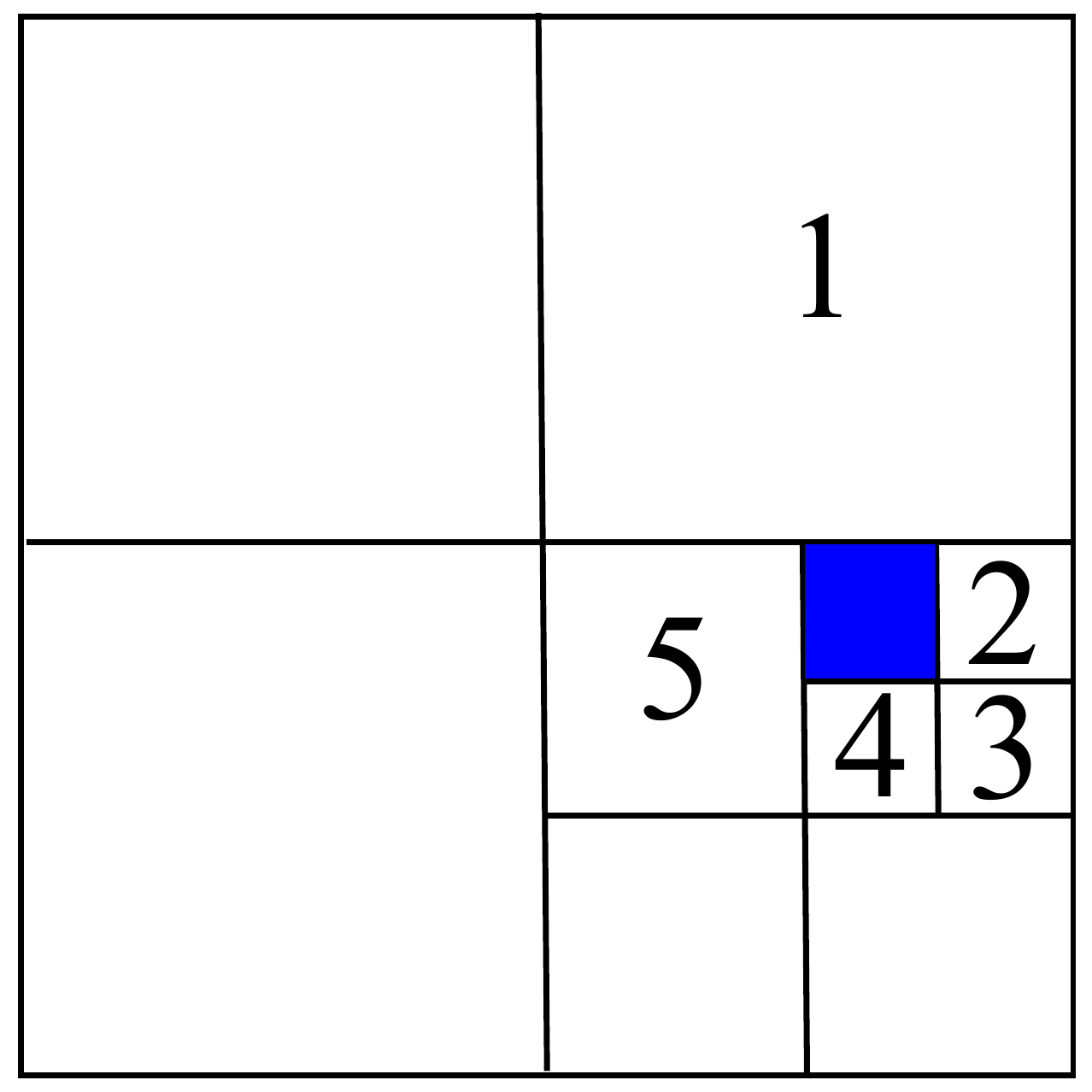}
    % \caption{ }
    % \label{fig:grid2d_worst_case}
  \end{minipage}
  \caption{ \footnotesize{Two particular configurations of the grid
      around a cell $\Omega_j$. Cell $\Omega_j$ is color-filled, while the neighbors are numbered. In these examples cell $\Omega_j$ has eight (left panel, uniform grid) or five (right panel, adaptive grid) neighbors, the latter being the case with the minimum number of neighbors.}  }
  \label{fig:grid2d}
\end{figure}

In two dimensions the extension of the reconstruction described
in~\cite{LPR:2001} to the case of adaptive grids is not
straightforward, so we describe it here in some detail. First, we define the set $\mathcal{N}_j$ of
neighbors of $\Omega_j$ as (see two examples in Fig.\ \ref{fig:grid2d}):
\begin{equation}\label{neigh}
 \mathcal{N}_j=\left\{ k \neq j \colon \partial \Omega_k \cap \partial \Omega_j \neq \emptyset \right\}.
\end{equation}
%In the examples in Fig.\ \ref{fig:grid2d} the color filled cell $\Omega_j$ has eight (left panel, uniform grid) or five (right panel, adaptive grid) neighbors.

The optimal polynomial $\Popt$ is chosen among the quadratic polynomials
matching the cell average $U_j$:
\begin{multline}\label{Popt2d}
 \Popt = U_j + p_x (x-x_j)+p_y(y-y_j) \\
+ \frac{1}{2} p_{xx} \left( (x-x_j)^2-\frac{h_j^2}{12} \right) +
\frac{1}{2} p_{yy} \left( (y-y_j)^2-\frac{h_j^2}{12} \right) + p_{xy} (x-x_j)(y-y_j).
\end{multline}
The coefficients $p_x$, $p_y$, $p_{xx}$, $p_{xy}$ and $p_{yy}$ are
determined imposing that $\Popt$ fits the cell averages of the cells in
$\mathcal{N}_j$ in a least-square sense. This can be realized if
$|\mathcal{N}_j| \geq 5$. Such condition is always satisfied in a quad
tree mesh. In fact, cell $\Omega_j$ has at least $3$ neighbors among
the children of its father cell and two edges on
the boundary of the father cell and thus at least two more neighbors
outside the father cell (see Fig.\ \ref{fig:grid2d}, in particular the
figure on the right for the case in which the cardinality of $\N_j$ is
as minimum as possible).
%By a similar argument, one can show that the
%neighborhood of a cell in a three-dimensional octree mesh contains at
%least $10$ cells ($2^d+d-1$ in the $d$-dimensional case), which
%provide enough data to determine the $9$ unknown coefficients
%($\left( \begin{array}{cc} 2+d \\ d \end{array} \right)-1$ in the
%$d$-dimensional case) of a second order polynomial in three space
%variables matching exactly the cell average $U_j$.
This result can be generalized in $d$ dimensions, where a cell has at least $2^d-1$ neighbors among its brother cells and shares at least $d$ neighbors with its father cell, making the minimum number of neighbors equal to $2^d+d-1$, while the unknown coefficients of a second degree polynomial in $d$ variables matching the cell average $U_j$ are $\left( \begin{array}{cc} 2+d \\ d \end{array} \right)-1$. It can be easily proved that $2^d+d-1 \geq \left( \begin{array}{cc} 2+d \\ d \end{array} \right)-1$ for all $d \geq 1$.

The equations to be satisfied in a least-square sense are:
\begin{equation}\label{eqLS}
\frac{1}{h_k^2} \int_{\Omega_k} \Popt\: dx = U_k, \qquad k \in \mathcal{N}_j.
\end{equation}
From \eqref{Popt2d} and \eqref{eqLS} we obtain the possibly over-determined
system $A\vec{c}=\vec{r}$, where $c=[p_x,p_y,p_{xx},p_{xy},p_{yy}]^T$,
$r_k=U_k-U_j$ denotes the component of $\vec{r}$ corresponding to neighbor $k$,
while $A$ is a $|\mathcal{N}_j|
\times 5$ matrix whose row corresponding to neighbor $k$ is:
\[
\left(
\begin{array}{c}
 x_k-x_j\\
y_k-y_j \\
\frac{1}{2} \left( (x_k-x_j)^2+\frac{1}{12}(h_k^2-h_j^2) \right) \\
(x_k-x_j) (y_k-y_j) \\
\frac{1}{2} \left( (y_k-y_j)^2+\frac{1}{12}(h_k^2-h_j^2) \right)
\end{array}
\right)^T
\]
We note in passing that the idea of employing least squares
fitting of the data in the same context has been exploited, in different
fashions, at least in \cite{Abgrall:1994,Furst:2006,LiRen:2012}.

Following the same procedure as in the one dimensional case, we
introduce four linear functions $P_\gamma$, $\gamma=1,\ldots,4$ in
such a way that $P_\gamma$ matches the cell average $U_j$ exactly 
\begin{equation}\label{PiUj}
 P_\gamma = U_j+p^\gamma_x (x-x_j)+p^\gamma_y(y-y_j)
\end{equation}
and the cell averages of a suitable subset of $\mathcal{N}_j$ in a
least-square sense. In practice, we choose these four subsets (stencils) along the direction of each of the four vertices of the central cell.
Let us start describing the choice of these stencils made in \cite{LPR:2001} for the uniform grid case. We introduce the following sets:

\begin{equation}\label{sets}
\begin{split}
 \mathcal{N}_j^E=\left\{ k \in \mathcal{N}_j \colon  x_k\geq x_j \right\}, \qquad &
\mathcal{N}_j^W=\left\{ k \in \mathcal{N}_j \colon  x_k\leq x_j \right\}, \\
 \mathcal{N}_j^N=\left\{ k \in \mathcal{N}_j \colon  y_k\geq y_j \right\}, \qquad &
\mathcal{N}_j^S=\left\{ k \in \mathcal{N}_j \colon  y_k\leq y_j \right\}.
\end{split}
\end{equation}
and define the four stencils:
\begin{equation}\label{stencils}
 \mathcal{N}_j^{\alpha \beta}=\mathcal{N}_j^{\alpha} \cap \mathcal{N}_j^{\beta}, \mbox{ with } \alpha \in \left\{ N,S\right\}, \; \beta \in \left\{ W,E\right\}.
\end{equation}
Referring to the left grid of Fig.\ \ref{fig:grid2d}, sets \eqref{sets} are:
\begin{equation*}
\begin{split}
 \mathcal{N}_j^E=\left\{ \Omega_k \colon k=1,2,3,4,5 \right\}, \qquad &
\mathcal{N}_j^W=\left\{ \Omega_k \colon k=1,5,6,7,8 \right\}, \\
 \mathcal{N}_j^N=\left\{ \Omega_k \colon k=1,2,3,7,8 \right\}, \qquad &
\mathcal{N}_j^S=\left\{ \Omega_k \colon k=3,4,5,6,7 \right\},
\end{split}
\end{equation*}
and stencils \eqref{stencils} are:
\begin{equation*}
\begin{split}
 \mathcal{N}_j^{NE}=\left\{ \Omega_k \colon k=1,2,3 \right\}, \qquad &
\mathcal{N}_j^{NW}=\left\{ \Omega_k \colon k=1,7,8 \right\}, \\
 \mathcal{N}_j^{SE}=\left\{ \Omega_k \colon k=3,4,5 \right\}, \qquad &
\mathcal{N}_j^{SW}=\left\{ \Omega_k \colon k=5,6,7 \right\}.
\end{split}
\end{equation*}
Observe that if we apply the same procedure in the case of adaptive grids, we could have some stencil with only one cell. In fact, referring to the adaptive structure of the right grid of Fig.\ \ref{fig:grid2d}, we would have:
\begin{equation*}
\begin{split}
 \mathcal{N}_j^E=\left\{ \Omega_k \colon k=2,3,4 \right\}, \: &
 \mathcal{N}_j^N=\left\{ \Omega_k \colon k=1,2 \right\} \: 
 \Longrightarrow
 \mathcal{N}_j^{NE}=\left\{ \Omega_k \colon k=2 \right\}.
\end{split}
\end{equation*}
Since having only one cell in a stencil is not enough to determine a linear function \eqref{PiUj}, we slightly modify the definition of the four sets \eqref{sets} (in such a way we include cell $\Omega_1$ in the $\mathcal{N}_j^E$ stencil of the right grid of Fig.\ \ref{fig:grid2d}). The new definition is:
\begin{eqnarray*}
 \mathcal{N}_j^E=\left\{ k \in \mathcal{N}_j \colon  x_k+ \displaystyle \frac{h_k}{2}\geq x_j \right\}, \qquad &
\mathcal{N}_j^W=\left\{ k \in \mathcal{N}_j \colon  x_k- \displaystyle \frac{h_k}{2}\leq x_j \right\}, \\
 \mathcal{N}_j^N=\left\{ k \in \mathcal{N}_j \colon  y_k+ \displaystyle \frac{h_k}{2}\geq y_j \right\}, \qquad &
\mathcal{N}_j^S=\left\{ k \in \mathcal{N}_j \colon  y_k- \displaystyle \frac{h_k}{2}\leq y_j \right\}.
\end{eqnarray*}
and the four stencils are defined as in \eqref{stencils}.
Now, referring to the adaptive structure of the right grid of Fig.\ \ref{fig:grid2d}, we have:
\begin{equation*}
\begin{split}
 \mathcal{N}_j^E=\left\{ \Omega_k \colon k=1,2,3,4 \right\}, \: &
 \mathcal{N}_j^N=\left\{ \Omega_k \colon k=1,2,5 \right\} \: 
 \Longrightarrow
 \mathcal{N}_j^{NE}=\left\{ \Omega_k \colon k=1,2 \right\}.
\end{split}
\end{equation*}

Let us rename the four stencils as $\mathcal{N}_j^{\gamma}$,
$\gamma=1,\ldots,4$.  The pairs of coefficients $p^\gamma_x$ and
$p^\gamma_y$ of \eqref{PiUj} are determined solving the system
$A^\gamma c^\gamma=r^\gamma$
in a least-square sense, with $c^\gamma=[p^\gamma_x,p^\gamma_y]^T$, $r^\gamma=[U_k-U_j, \: k \in
\mathcal{N}_j^\gamma]$, and $A^\gamma$ is a $|\mathcal{N}_j^\gamma | \times
2$ matrix whose $k-$th row is:
\[
\left(
\begin{array}{c}
 x_k-x_j\\
y_k-y_j
\end{array}
\right)^T
\]
We observe that also in this case the least-square problem is
well-determined because $|\N_j^\gamma | \geq 2$, since we have at
least one cell on the other side of each cell's edge (see Fig.\ \ref{fig:grid2d}).
\begin{remark}
We observe that in the uniform grid case the scheme does not reduce to the 2D CWENO described in \cite{LPR:2001}. In fact here each plane is obtained solving an over-determined system of 
four equations (in an adaptive grid $\left|\mathcal{N}_j^{\gamma}\right|+1$ equations, but $\left|\mathcal{N}_j^{\gamma}\right|=3$ in the uniform case) and three coefficients, while in \cite{LPR:2001} each plane was obtained by imposing only three conditions.  
\end{remark}
The parabola $P_0$ is determined by the relation:
\begin{equation}\label{P0}
 \Popt = \alpha_0 P_0 + \sum_{i=1}^4 \alpha_\gamma P_\gamma.
\end{equation}
By the same argument of the one dimensional case, the choice of
coefficients $\alpha_\gamma$ is arbitrary.  In practice we use
$\alpha_0=1/2$ and $\alpha_\gamma=1/8$, $\gamma=1,\ldots,4$.  The
final reconstruction is:
\begin{equation}\label{finalrec}
 P= \tilde{\alpha}_0 P_0 + \sum_{\gamma=1}^4 \tilde{\alpha}_i P_\gamma
\end{equation}
with
\[
 \tilde{\alpha}_\gamma = \frac{\omega_\gamma}{\sum_{\delta=0}^4 \omega_\delta}, \qquad \omega_\gamma = \frac{\alpha_\gamma}{(\epsilon+\IS_\gamma)^2}, \quad \gamma=0,\ldots,4.
\]
The smoothness indicators, $\IS_\gamma$, are~\cite{LPR:2001}:
\begin{equation}\label{IS2d}
\IS_\gamma=\sum_{|\alpha|=1}^2 \int_{\Omega_j} h_j^{2(|\alpha|-1)} (P_\gamma^{\alpha})^2 d \Omega, \qquad \gamma=0,\ldots,4.
\end{equation}
where $\alpha=(\alpha_x,\alpha_y)$ is a multi-index denoting the derivatives.
Let us rewrite the polynomials $P_\gamma$, $\gamma=0,\ldots,4$ as:
\begin{eqnarray*}%\label{Pi2d}
 P_{\gamma} &=& U_j + p_x^\gamma (x-x_j)+p_y^\gamma(y-y_j) + \frac{1}{2} p_{xx}^\gamma \left( (x-x_j)^2-\frac{h_j^2}{12} \right)\\
 &+& \frac{1}{2} p_{yy}^\gamma \left( (y-y_j)^2-\frac{h_j^2}{12} \right) + p_{xy}^\gamma (x-x_j)(y-y_j).
\end{eqnarray*}
A direct computation of \eqref{IS2d} yields:
\[
\IS_\gamma=h_j^2(p_x^\gamma)^2+h_j^2(p_y^\gamma)^2+\frac{13}{12} h_j^4(p_{xx}^\gamma)^2+\frac{7}{6} h_j^4(p_{xy}^\gamma)^2+\frac{13}{12} h_j^4(p_{yy}^\gamma)^2.
\]
Observe that $\pxx^\gamma=\pxy^\gamma=\pyy^\gamma=0$ for $\gamma=1,\ldots,4$.

\subsection{Treatment of the boundaries}
Boundary conditions are treated as follows. We create one layer of ghost cells around the computational domain. The size of each ghost cell matches the size of the adjacent cell in the domain $\Omega$. The value of the cell average of the field variables in the ghost cell is determined by the boundary conditions (e.g.\ free flow BC, reflecting BC, and so on). Because of the compactness of the scheme, only one layer is necessary in order to perform the reconstruction in the physical cell near the boundary. The field variable on the outer side of the physical boundary, necessary to compute the numerical flux, is computed from the field variable on the inner side by applying the boundary conditions again. For example, in free flow boundary conditions the inner and outer values of the field variables are the same, while for reflecting BC in gas dynamics, density, pressure and tangential velocity on both side are equal, while the outer value of the normal velocity is the opposite of the inner one.

\section{Fully discrete scheme and adaptivity}
In this section we describe the fully discrete scheme that we use to test the CWENO reconstruction, that is we specify the time integration procedure and the strategy for adaptive mesh refinement. 
First we introduce some notation for semidiscrete schemes.

\begin{figure}
\begin{center}
\includegraphics[width=0.5\textwidth]{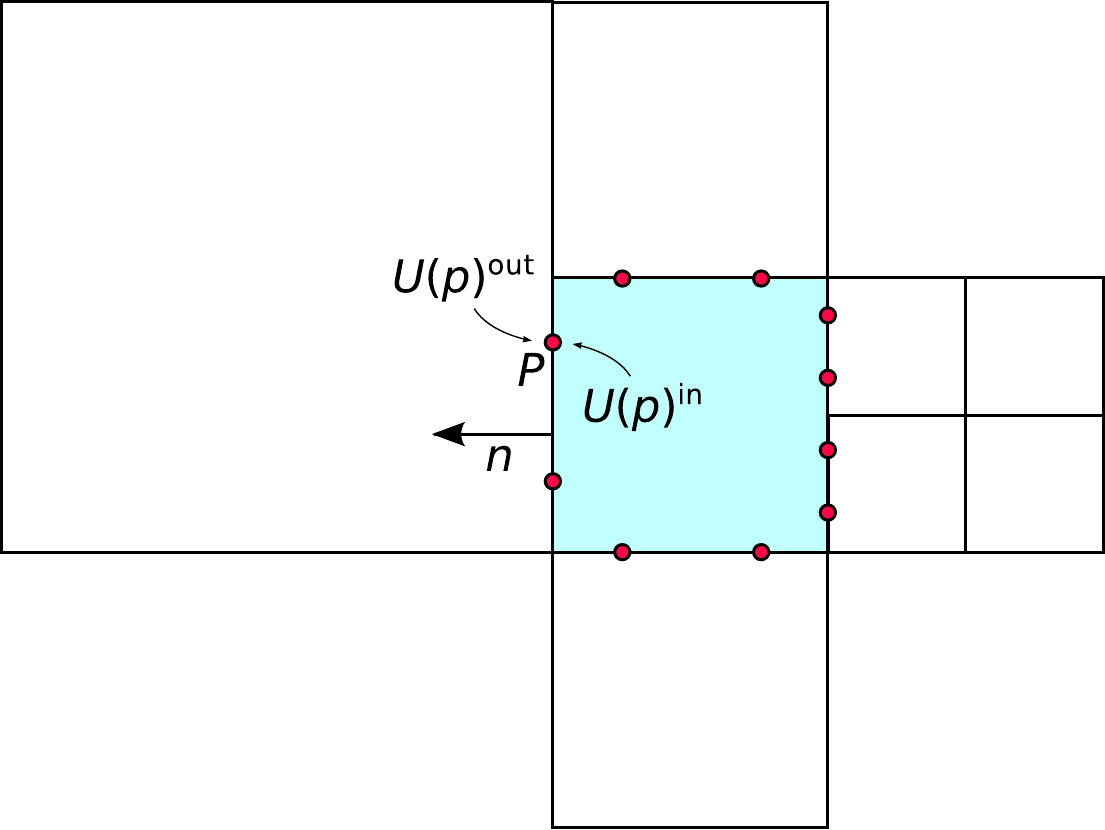}
\end{center}
\caption{Location of quadrature points for the computation of the numerical fluxes. The colored cell has five adjacent neighbors, and therefore uses the fluxes computed in the ten quadrature nodes marked by small circles.}
\label{fig:gausspoints}
\end{figure}

Denote with $U_j^n$ the approximate cell average of the solution at
time $t^n$ in cell $\Omega_j$. Then a first order (in time) scheme for system \eqref{eq:finvol} may
be written as
\begin{equation*}
U^{n+1}_j = U^n_j - 
\frac{\DT}{|\Omega_j|}
\mathcal{Q}(\partial\Omega_j;\mathcal{F})
\end{equation*}
where $\mathcal{Q}$ is a suitable quadrature formula and $\mathcal{F}$
are the numerical fluxes. 
Let ${\mathcal N}^A_j$ denote the set of {\em adjacent} neighbors of cell $j$ (cells that share a segment with $\Omega_j$). Then we split the boundary $\partial \Omega_j$ as 
\[
   \partial \Omega_j = \bigcup_{i\in\mathcal{N}^A_j} \partial \Omega_j \cap \partial \Omega_i.
\]
We choose a Gauss quadrature
on each segment of this decomposition and compute the numerical
fluxes at each quadrature node as
$\mathcal{F}(f\cdot\vec{n};U(x)^{\text{in}},U(x)^{\text{out}})$, where
$f:\R^d\to\R^m$ is the exact flux function, $U(x)^{\text{in}}$ is the
reconstruction in cell $j$ evaluated at $x$ and $U(x)^{\text{out}}$ is
obtained evaluating at $x$ the reconstruction in the nearby cell. Each numerical flux computed on a quadrature node is used on the two adjacent cells it belongs, thus automatically guaranteeing conservation. Note that
the concept of nearby cell is unambiguous if the quadrature nodes are
not on the vertices of the cell, as is the case for Gauss-type
formulas. An example is depicted in Fig.\ \ref{fig:gausspoints}.

An important issue in the construction of an efficient adaptive method 
concerns the time stepping: because of hyperboic CFL restriction, time step 
has to be of the order of the local mesh size, therefore one would have $\DT = O(H)$ in smooth regions, and $\DT = O(h)$ near singular regions. Time advancement procedures that can employ different time step length in different regions of the computational domain are called local timestepping schemes.
Originally local time stepping methods have been developed on uniform grids, for problems with highly non uniform propagation speeds, in order to guarantee an optimal CFL condition on the whole domain (see \cite{OS:1983}). Later they have been adopted to improve the efficiency of schemes based on non uniform grids. 
Different strategies have been adopted to implement local time stepping. In CLAWPACK, for example, the computational domain is discretized with a coarse uniform grid of mesh size $H$, which contains patches of more refined grid where necessary. At macroscopic time step $\DT=O(H)$ the solution is updated on the coarse grid, and refined on the finer grid and on the coarse cells adjacent to the fine grid region, making sure that conservation is guaranteed (see \cite{BergerLeVeque98} for details). More direct second order AMR methods with adaptive time step have been adopted, see for example \cite{Kir:2003,CS:2007}, where a detailed analysis is performed of second order local time step methods, \cite{PS:2011:entropy} for an application of local timestepping in a very similar setting, or \cite{LGM:2007}, where an asyncronous time step strategy is adopted, in which the cell with the smallest time is the one that is advanced first. 

In the present paper we shall not discuss local time step, since the main point of the paper is to extend and analyze Compact WENO discretization to adaptive grids. Local time stepping for third order schemes is non-trivial, and will be the subject of a future paper. The time step $\DT$ will thus be chosen in order to satisfy the CFL
condition everywhere in the domain of the PDE.
% This usually implies that
%small cells dictate the timestep, unless multirate timestepping is
%employed (see e.g. \cite{Kir:2003,CS:2007,PS:2011:entropy} for second order examples or \cite{LGM:2007} for a more general approach). For the sake of simplicity,
%in this paper we will restrict ourselves to classic explicit
%TVD-RungeKutta schemes as the focus of the present paper is on the
%CWENO reconstruction.
Thus our fully discrete scheme will be written as
\begin{equation}\label{eq:RK}
U^{n+1}_j = U^n_j - 
\frac{\DT}{|\Omega_j|}
\sum_{i=0}^{\sigma}
b_i \mathcal{Q}(\partial\Omega_j;\mathcal{F}^{(i)})
\end{equation}
where the stage fluxes $\mathcal{F}^{(i)}$ are computed by applying
the same formula as above to the stage values of the explicit Runge-Kutta
scheme
\begin{equation*}\label{eq:RKstage}
 U^{(i)}_j = U^n_j - 
\frac{\DT}{|\Omega_j|}
\sum_{k=0}^{i-1}
a_{ik} \mathcal{Q}(\partial\Omega_j;\mathcal{F}^{(k)})
\end{equation*}
Here $(a_{ij},b_i)$ denote the coefficient of the Butcher tableaux of
a Strong Stability Preserving Runge-Kutta scheme with $\sigma$ stages
(see \cite{GKS:book:SSP}). In all the tests of this paper we employ
the Local Lax-Friedrichs numerical fluxes.

\subsection{Error estimators/indicators}
Several techniques can be adopted to decide where refine or derefine locally the mesh. Most of them are based on the use of local error indicators, such as, for example, discrete gradients and discrete curvature \cite{ArvDelis:2006}, interpolation error \cite{BergerLeVeque98}, residuals of the numerical solution \cite{KK:2005} or of the entropy \cite{Ohlberger:2009,Puppo:2003:entropy}.
In the present paper we employ the numerical entropy production, that was introduced in \cite{Puppo:2003:entropy} for central schemes and later extended in \cite{PS:2011:entropy} to unstaggered finite volume schemes of arbitrary order. 
The motivation for this choice is that the numerical entropy production is naturally available for any system of conservation laws with an entropy inequality, it scales as the truncation error in the regular regions, and its behavior allows to distinguish between contact discontinuities and shocks.

In order to construct the indicator, one considers an entropy pair
$(\eta,\psi)$ and chooses a numerical entropy flux $\Psi$ compatible
with the exact entropy flux $\psi$ in the usual sense that
$\Psi(u,u)=\psi(u)$ and $\Psi$ is at least Lipshitz-continuous in each
entry. Then one forms the quantity
\begin{equation}\label{eq:entprod}
  S^n_j = \frac1{\DT} 
  \left[ \langle\eta(U^{n+1})\rangle_j - \langle\eta(U^{n})\rangle_j 
    + \frac{\DT}{|\Omega_j|}
    \sum_{i=0}^{\sigma}
    b_i \mathcal{Q}(\partial\Omega_j;\Psi^{(i)})
  \right]
\end{equation}
where we have denoted with $\langle\cdot\rangle_j$ the operation of
averaging on the cell $\Omega_j$ and $\Psi^{(i)}$ denote the numerical
entropy fluxes computed using the reconstruction of the $i$-th stage
value.

In one spatial dimension, \cite{PS:2011:entropy} show that, if the
solution is locally smooth, $S_j^n=O(h^r)$ with $r$ equal to the
minimum between the order of the scheme and the order of the
quadrature formulas used to compute the cell averages of the entropy;
on the other hand, $S_j^n=O(h)$ (resp. $O(1/h)$) if there is a contact
discontinuity (resp. shock) in $\Omega_j$. Of course, for a second
order scheme it is enough to employ the midpoint rule
$\langle\eta(u^{n})\rangle_j=\eta(U^n_j)$ as in
\cite{PS:2011:entropy}, whereas in order to observe a third order
scaling of $S_j^n$ as $h\to0$, one has to employ a quadrature formula
of higher order and thus evaluate the reconstructions of $u^n$ and
$u^{n+1}$ at the quadrature points. Since the CWENO reconstruction
yields a polynomial with uniform accuracy in the whole $\Omega_j$,
this entails simply an evaluation of the already computed polynomial
and does not involve other reconstruction steps or extra sets of
weights. Of course the reconstruction of $u^{n}$ is already available
and the reconstruction of $u^{n+1}$ would be computed in any case at
the beginning of the next time step.

\begin{figure}
\begin{center}
\includegraphics[width=0.99\textwidth]{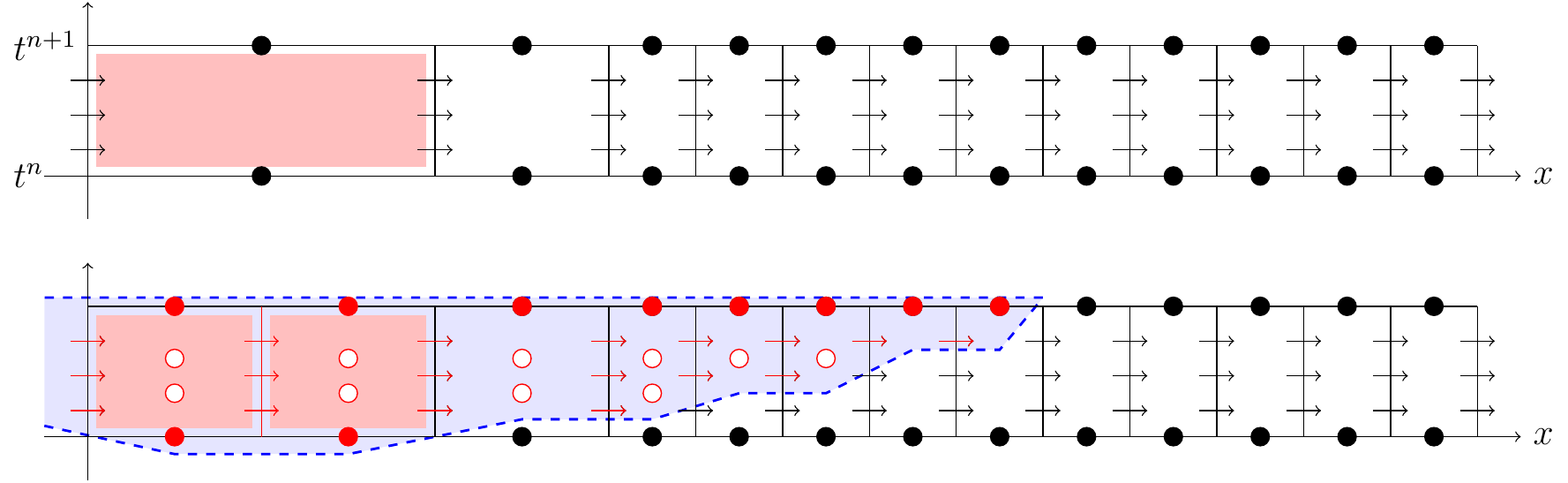}
\end{center}
\caption{Recomputation after refinement. The top row depicts the tentative computation of the time advancement after which the shaded cell is marked for refinement. Circles represents the cell averages and the arrows the numerical fluxes computed during the three Runge-Kutta stages. The bottom row depicts the recomputation performed after splitting the shaded cell. The dashed line indicates the numerical domain of dependence originating from the split cell. The stage values of the Runge-Kutta scheme that are recomputed are indicated by empty circles. The other stage values are not indicated for easier reading.}
\label{fig:recompute}
\end{figure}

Following ideas from \cite{PS:2011:entropy} and extending them to
the case of $d$ space dimensions, we construct an adaptive mesh refinement
scheme as follows:
\begin{itemize}
\item we start from a uniform coarse grid with $N_0$ cells, called the level-0 grid.
\item any single cell of the grid may be replaced by $2^d$ equal cells
  and this operation may be
  performed recursively, obtaining a computational grid that is
  conveniently stored in a binary, quad- or oct-tree depending on $d$. Note that these grids will have in
  general hanging nodes, but this poses no problem to finite-volume
  discretizations. Let's call {\em level} of a cell its depth in the
  tree; obviously a cell of level $l$ has diameter $2^l$ times smaller than its level-0 ancestor.
\item at the end of each timestep, the quantity $S^n_j$ is computed in
  every cell. If
  it is bigger than a threshold $\Sref$ and if the level of refinement of $\Omega_j$ does
  not equal the maximum refinement level allowed in the grid, the cell
  is refined. The cell averages in the newly created cells are set
  by averaging the reconstruction of $\{U^n_{\cdot}\}$\footnote{As an exception, during the
    first time step, if the initial condition is known analytically,
    it is more accurate to use the analytic expression to set the cell
    averages in the newly created cells.} and the timestep recomputed
  locally, i.e. $U^{n+1}_i$ 
  is recomputed only in the numerical domain of dependence of the cell $\Omega_j$ (see Fig.\ \ref{fig:recompute}).
\item the solution $\{U^{n+1}_j\}$ is accepted when no refinement is
  required nor possible by the conditions on the size of $S^n_j$ and on
  the level of the cells. Note that since $S^n_j$ diverges on shocks when the grid is refined,
  fixing a maximum refinement level is necessary.
\item a coarsening pass checks if all $2^d$ direct children of a
  previously refined cell have an entropy production lower than a given coarsening threshold, i.e.\
  $S^n_j<\Scoa$ and if so,
  replaces the $2^d$ children with their ancestor cell, where it sets
  $U^{n+1}_j$ equal to the average of the cell averages in the
  children. As in~\cite{PS:2011:entropy} we employ $\Scoa=\Sref/2^{p+1}$. At this point a new timestep starts.
\end{itemize}

Our code makes use of the DUNE interface~\cite{dunegridpaperI:08,dunegridpaperII:08,dune-web-page} 
to achieve grid-independent and dimension-independent coding of the
numerical scheme in C++. Such interface is able to adopt several kinds of grid packages, including the ALUGRID library~\cite{ALUGRID},
which is the one that has been adopted in the two dimensional simulations. 
More precisely, a {\em dune-module} called
{\tt dune-fv} was written by the first author to provide generic
interfaces to explicit Runge-Kutta, numerical fluxes, adaptive
strategy, reconstructions, as well as the implementation of the
classical schemes; the CWENO reconstruction was incorporated in {\tt dune-fv}
by the second author. The source code, licensed under GPL terms, allows
to build adaptive mesh refinement schemes in one or two space
dimensions, with spatial and temporal order up to three. Most
components are easily interchangeable and the source allows for easy
experimentation with different combinations of time-stepper,
reconstructions, numerical fluxes, error indicator, etc
\cite{dunefv}

\section{Order of accuracy}
In this section we provide a scaling argument in support of the use  of the third order scheme. We distinguish between the computation of regular solutions and piecewise smooth solutions. In the case of regular solutions, the observed order  of accuracy of the method should be equal  to the theoretical one, giving an error that scales like $N^{-r/d}$, where $N$ is the total number of cells, $r$ is the (space and temporal) order of the scheme and $d$ the number of space dimensions.

In this case, h-adaptivity helps by automatically refining in the regions of smaller space scales (see e.g. the tests of Fig.\ \ref{fig:lintra2} for one space dimension and Fig.\ \ref{fig:swirl} in two space dimensions). The ratio between the largest and the smallest mesh size depends on the ratio between the largest and the smallest macroscopic scale of the physical system. The CFL condition implies that the time discretization is bounded by the space mesh scale, which suggest that the optimal choice is obtained by using the same order of accuracy in both space and time.

The situation is very different when discontinuities are present. Let us discuss separately the cases of one and of higher spatial dimensions. 

\subsection{One dimensional case}
Let us assume that we want to solve a problem whose solution is piecewise smooth. For simplicity, we assume that the spatial and temporal order of accuracy are the same.
More precisely, let us assume that we want to solve a generalized Riemann problem:

\begin{eqnarray*}
&& u_t + f(u)_x = 0, \\ 
&&u(x,0)=u_0(x)+H(x)\delta u(x),
\end{eqnarray*}
where $u_0$ and $\delta u$ are smooth functions and $H$ is the Heaviside function 
\[
H(x)=
\left\{
\begin{array}{cc}
1 & \mbox{ if } x \geq 0 \\
0 & \mbox{ if } x < 0 
\end{array}.
\right.
\]
For some time the solution of such problem will consist of a piecewise smooth function, in which singularities are localized in a few points. For example, for Euler equations in gas dynamics we have a jump in all quantities at the shock, a jump in density at contact discontinuity and a jump in the first space derivatives at the end of the rarefaction region induced by the initial discontinuity.

Let $r$ denote the order of the scheme employed, $H$ the mesh size of the grid that we would use on a smooth region, $h$ the smallest mesh that we allow by refinement. During the evolution, the adaptive algorithm will refine around the singularity points. We can identify two regions: a smooth region, bounded away from singularities, and a singular region, whose size is $O(H)$ around the singularities. Assume that a non adaptive scheme resolves each singularity in $\nu$ grid points.\footnote{For simplicity we assume $\nu$ is independent on the nature of the singularity, which of course is not true in general.} Then the size of the singular region is $\mathcal{L}_s=\nu H$. Near the singularity the method reduces to first order. If we want that the lack of accuracy caused by the low order scheme does not pollute the solution in the smooth region, we have to require that $h=O(H^r)$. The motivation for this choice is the following.
The solution in the regular region is affected by the presence of the singularity, with an error which is the maximum between the truncation error of size $O(H^r)$ due to the method and the error induced by the first order treatment of the singularity, which is $O(h)$. Therefore the error, even in smooth regions, is $\max(O(H^r),O(h))$, which is the motivation for suggesting the optimal choice $O(h) = O(H^r)$.

%In addition, we assume hyperbolic CFL condition, i.e.\ $\Delta t = O(H)$ in the smooth region, and use a time integration of order $p$ as well.
In this way the overall accuracy of the solution in the regular regions should not be affected by the presence of the singularity, in the sense that the error introduced by the singularity is of the same order of the local truncation error and one expects that the error after a finite time is $O(H^r)$. This is true in local timestepping schemes where $\Delta t_{\max} = O(H)$ and even more so in ths present case where we employ global time stepping, i.e. $\Delta t = O(h)$.

\begin{figure}
\begin{minipage}[t]{0.55\textwidth}
\vspace{0pt}
\includegraphics[width=\textwidth]{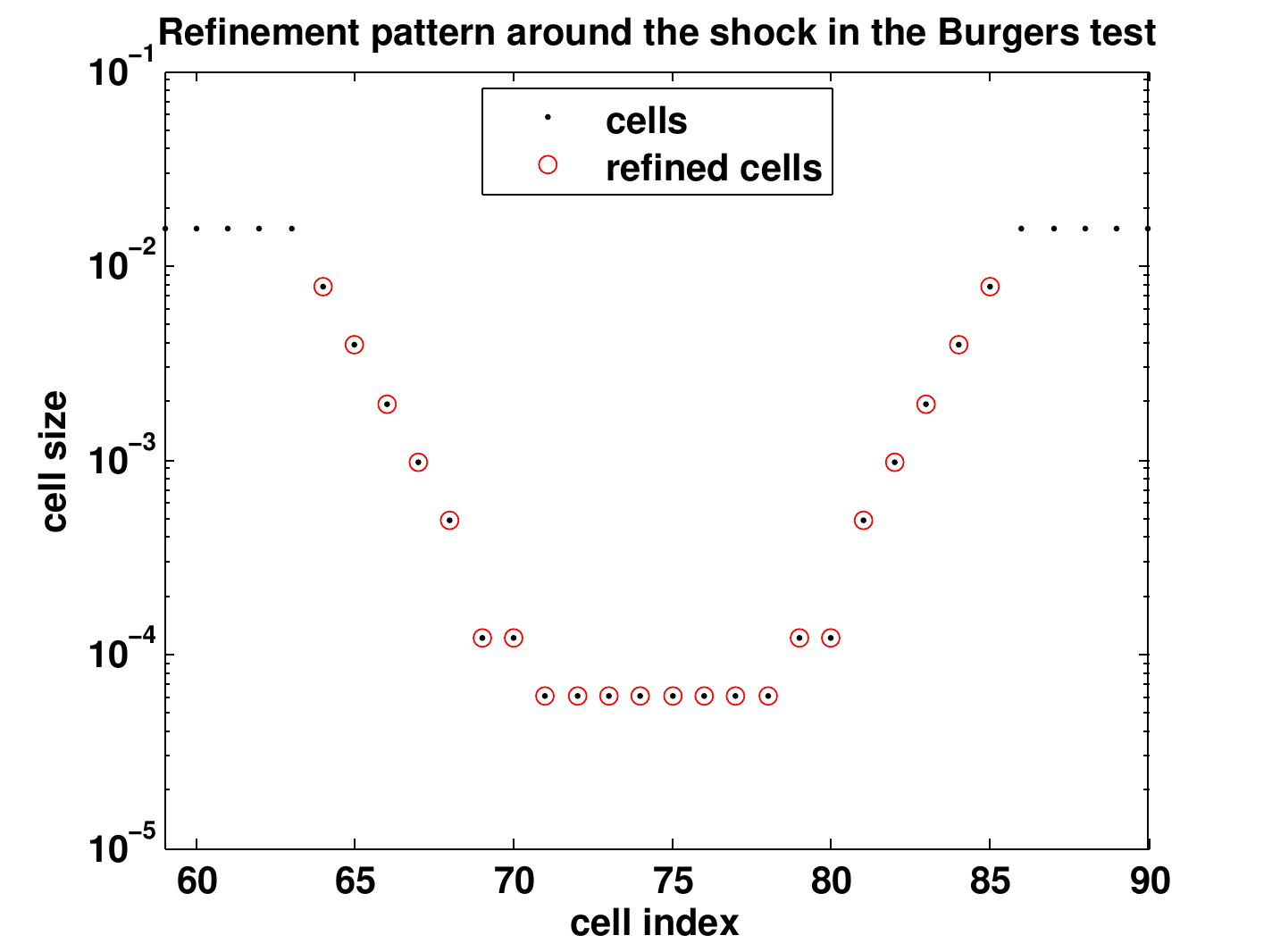}
\end{minipage}
\hfill
\begin{minipage}[t]{0.35\textwidth}
\vspace{26pt}
\begin{tabular}{|c|c|} \hline
Number of & Number of\\
levels & refined cells \\ \hline
3 & 6  \\
5 & 10 \\
7 & 16 \\
9 & 22 \\ \hline   
\end{tabular}
\end{minipage}
\caption{Refined cells around the shock at $x=0$ for the Burgers test with standing shock (at the final time of the simulation). 
Cell size vs index when using $L=9$ refinement levels (left); number of cells involved in refinement vs refinement levels (right).} 
\label{fig:shockcells}
\end{figure}

Let us check the dependence of the error on the total number of cells. Let $\mathcal{L}$ be the length of the computational domain. Then the number of cells in the large region is $N_R=\displaystyle \frac{\mathcal{L}-N_s \nu H}{H} = \displaystyle \frac{\mathcal{L}}{H}-N_s \nu$, where $N_s$ denotes the number of singular points.
%Without loss of generality we can assume $\mathcal{L}=1$.
The number of cells in the singular region depends on the number of refinement levels 
$L = \log_2 \displaystyle \frac{H}{h}$. In fact, if we assume that we have a shock located at a given point, we keep refining the cells until we reach the minimum cell size. From the numerical experiments (see Fig.\ \ref{fig:shockcells} for the case of the Burgers equation after shock formation) it appears that the number of refined cells around one shock is approximately equal to $L$ plus a number of smallest cells located at the shock. Such a number appearsh to be mildly increasing with $L$, so that the total number of points cells refined due to the presence of a shock is much smaller than the total number of cells.
Since $h=c H^r + \text{h.o.t.}$ for some constant $c$, we have:
\[
L\simeq\log_2 \frac{H}{cH^r} = \log_2 \frac{1}{c} - (r-1) \log_2 H.
\]
The number of cells in the singular region $N_I$ is therefore negligible with respect to the number of cells in the smooth region, $N_I = O(L) \ll N_R$. Therefore, the total number of cells $N=N_R+N_I \approx \mathcal{L}/H$. As a consequence, the loglog plot of Errors vs.\ $N$ should show the classical slope close to $-r$. Notice, however, the slight dependence of the number of levels on $H$. In practice, if we change $H$ with $H/2$, in order to obtain an error decay rate of $r$, we also have to increase the number of levels by an amount $r-1$. This is because the error in the smooth region becomes $2^r$ times smaller, so $h$ has to be divided by $2^r$ (Table \ref{table:nooflevels1D}).
\begin{table}
\centering
\begin{tabular}{|c|c|c|}
largest $\Delta X$ & smallest $\Delta x$ & no. of levels \\
\hline
$H_1 = H$ & $h_1=h$ & $L_1= \log_2 (H/h)$ \\ %\displaystyle \frac{H}{h}$ \\
$H_2 = H/2$ & $h_2=h/2^r$ & $L_2= \log_2 (2^{r-1} H/h) = (r-1) + L_1 $
\end{tabular}
\caption{How the maximum number of levels has to change when we increase the resolution on the coarse grid by a factor $2$.}
\label{table:nooflevels1D}
\end{table}

\subsection{Higher dimensional case}
The situation is different in more space dimensions. Let us first analyze in detail the two dimensional case, which is the one relevant to the numerical tests of this paper. Let us assume that the singularities are concentrated on a line that evolves (see e.g. Fig.\ \ref{fig:riemann}).
Let $N$ denote the number of cells per direction on the coarsest grid. Without refinement there would be a total of $N_{\text{TOT}} = N^2$ cells. Assume now that we adopt a grid refinement until the finest cell becomes of size $h$. The cell size on the coarsest grid is $H=\mathcal{L}/N$. Since near the discontinuity the accuracy degrades to first order, it is natural to choose $h=c H^r$ for some constant $c$. The number of cells of size $h^2$ is proportional to the length of the discontinuity line divided by $h$, so $N_I = O(1/h)$, while $N_R= O(1/H^2)$. The total number of cells,  is thus $N_{\text{TOT}}\approx C_I/h + C_R/H^2$, for some constants $C_I$ and $C_R$.  Hence $N_{\text{TOT}} \approx C_I/(cH^p) + C_R /H^2$ and we have the following scaling:
\begin{description}
\item[$r<2$:] most cells lie in the regular region;
\item[$r=2$:] the number of cells in the regular and in the singular region have the same scaling;
\item[$r>2$:] most cells are near the singular line.
\end{description}
This argument suggests that, in the presence of discontinuities, it is not possible to observe an asymptotic behavior of the error better than $O(1/N_{\text{TOT}})$ even using higher order reconstruction in the smooth region. However, although the asymptotic behavior is the same, the use of a higher order reconstruction may produce, as we shall see, a considerably smaller error for the same coarsest mesh, or, alternatively, the same error can be obtained with a considerably smaller number of cells. That this is the case in practical tests can be appreciated in Fig.\ \ref{fig:swirl} for a smooth test and in Figures  \ref{fig:riemann} and \ref{fig:shockbubble} in the presence of shocks.

In three dimensions the singularities are in general concentrated on a two dimensional manifold. Therefore, even with a second order accurate method most of the cells will be of size $h$, since
 $N_{\text{TOT}}  \approx  \displaystyle \frac{C_I}{h^2} + \frac{C_R}{H^3}
 =\frac{C_I}{cH^{2r}} + \frac{C_R}{H^3}$. Therefore if $r=1$ most cells will be in the regular region, while with $r=2$ the majority of cells will be near the singular surface.

\section{Numerical tests}
In the first set of tests we want to assess the convergence order of
the CWENO reconstruction and of the fully discrete numerical scheme on
uniform and non-uniform grids.
We focus mainly on the choice of the parameter $\epsilon$ and on its dependence on the local grid size $h$. At variance with the case of a uniform grid, in which one can use a constant value of $\epsilon$, here it is crucial to incorporate its dependence on the local grid size, which may vary by orders of magnitude.

\begin{table}
\begin{minipage}[t]{0.49\textwidth}
Smooth test, $\epsilon=h_j$\\
\footnotesize
\begin{tabular}{r|cc|cc|}
N & $\|E\|_1$ & rate & $\|E\|_\infty$ & rate \\
\hline
  20 &   5.63e-03  &       & 2.17e-02 & \\
  40 &   5.04e-04  &  3.48 & 1.49e-03  &  3.87 \\
  80 &   4.97e-05  &  3.34 & 1.20e-04  &  3.64 \\
 160 &   5.40e-06  &  3.20 & 1.32e-05  &  3.18 \\
 320 &   6.22e-07  &  3.12 & 1.65e-06  &  3.00 \\
 640 &   7.43e-08  &  3.07 & 2.06e-07  &  3.00 \\
1280 &   9.07e-09  &  3.03 & 2.57e-08  &  3.00 \\
2560 &   1.12e-09  &  3.02 & 3.22e-09  &  3.00 \\
\end{tabular}
\end{minipage}
\begin{minipage}[t]{0.49\textwidth}
Smooth test, $\epsilon=10^{-6}$\\\
\footnotesize
\begin{tabular}{r|cc|cc|}
N & $\|E\|_1$ & rate & $\|E\|_\infty$ & rate \\
\hline
  20 & 1.30e-02  &        & 4.17e-02  &  \\
  40 & 2.26e-03  &  2.53  & 1.02e-02  &  2.03\\
  80 & 3.55e-04  &  2.67  & 2.52e-03  &  2.02\\
 160 & 4.95e-05  &  2.84  & 5.89e-04  &  2.10\\
 320 & 4.96e-06  &  3.32  & 5.44e-05  &  3.44\\
 640 & 4.36e-07  &  3.51  & 2.25e-06  &  4.60\\
1280 & 3.91e-08  &  3.48  & 1.31e-07  &  4.10\\
2560 & 3.27e-09  &  3.58  & 8.95e-09  &  3.87\\
\end{tabular}
\end{minipage}
\caption{Reconstruction errors for the smooth function $u_s$ on a uniform grid.}
\label{tab:eps:smooth}
\end{table}

%\begin{table}
%\begin{minipage}[t]{0.49\textwidth}
%Discontinuous test, $\epsilon=h_j$\\
%\footnotesize
%\begin{tabular}{r|cc|cc|}
%N & $\|E\|_1$ & rate & $\|E\|_\infty$ & rate \\
%\hline
%  20 & 2.14e-03  &      & 2.49e-02  &      \\
%  40 & 8.99e-04  & 1.25 & 2.02e-02  & 0.31 \\
%  80 & 3.16e-04  & 1.51 & 1.37e-02  & 0.56 \\
% 160 & 8.62e-05  & 1.87 & 7.28e-03  & 0.91 \\
% 320 & 1.78e-05  & 2.27 & 2.97e-03  & 1.29 \\
% 640 & 2.96e-06  & 2.59 & 9.82e-04  & 1.60 \\
%1280 & 4.30e-07  & 2.79 & 2.84e-04  & 1.79 \\
%2560 & 5.79e-08  & 2.89 & 7.64e-05  & 1.89 \\
%\end{tabular}
%\end{minipage}
%\begin{minipage}[t]{0.49\textwidth}
%Discontinuous test, $\epsilon=10^{-6}$\\\
%\footnotesize
%\begin{tabular}{r|cc|cc|}
%N & $\|E\|_1$ & rate & $\|E\|_\infty$ & rate \\
%\hline
%  20 & 4.22e-04  &      & 1.66e-03  &      \\
%  40 & 6.11e-05  & 2.79 & 4.16e-04  & 1.99 \\
%  80 & 6.54e-06  & 3.22 & 1.04e-04  & 2.00 \\
% 160 & 6.33e-07  & 3.37 & 2.60e-05  & 2.00 \\
% 320 & 6.09e-08  & 3.38 & 6.51e-06  & 2.00 \\
% 640 & 5.94e-09  & 3.36 & 1.63e-06  & 2.00 \\
%1280 & 6.30e-10  & 3.24 & 4.06e-07  & 2.00 \\
%2560 & 7.37e-11  & 3.10 & 1.01e-07  & 2.01 \\
%\end{tabular}
%\end{minipage}
%\caption{Reconstruction errors for a discontinuous function on a
%  uniform grid. The discontinuity is on a grid interface.}
%\label{tab:eps:jump}
%\end{table}

\begin{table}
\begin{minipage}[b]{0.49\textwidth}
Discontinuous test, $\epsilon=h_j$\\
\footnotesize
\begin{tabular}{r|cc|cc|}
N & $\|E\|_1$ & rate & $\|E\|_\infty$ & rate \\
\hline
   20 & 5.14e-03  &      & 7.51e-02 &      \\ 
   40 & 2.70e-03  & 0.93 & 7.98e-02 & -0.09 \\ 
   80 & 9.69e-04  & 1.48 & 4.74e-02 & 0.75 \\ 
  160 & 4.51e-04  & 1.11 & 4.99e-02 & -0.08 \\ 
  320 & 2.08e-04  & 1.11 & 5.28e-02 & -0.08 \\
  \hline 
  640 & 4.37e-06  & 5.57 & 9.82e-04 & 5.75 \\ 
 1280 & 6.36e-07  & 2.78 & 2.84e-04 & 1.79 \\ 
 2560 & 8.59e-08  & 2.89 & 7.64e-05 & 1.89 \\ 
\end{tabular}
\end{minipage}
\begin{minipage}[b]{0.49\textwidth}
Discontinuous test, $\epsilon=10^{-6}$\\\
\footnotesize
\begin{tabular}{r|cc|cc|}
N & $\|E\|_1$ & rate & $\|E\|_\infty$ & rate \\
\hline
   20 & 5.42e-03  &      & 1.02e-01 &      \\ 
   40 & 2.56e-03  & 1.08 & 1.00e-01 & 0.02 \\ 
   80 & 7.29e-04  & 1.81 & 5.78e-02 & 0.80 \\ 
  160 & 3.61e-04  & 1.01 & 5.77e-02 & 0.00 \\ 
  320 & 1.80e-04  & 1.00 & 5.77e-02 & 0.00 \\ 
  \hline
  640 & 6.42e-09  & 14.78 & 1.63e-06 & 15.11 \\ 
 1280 & 6.91e-10  & 3.22 & 4.06e-07 & 2.00 \\ 
 2560 & 8.12e-11  & 3.09 & 1.01e-07 & 2.01 \\ 
\end{tabular}
\end{minipage}
\caption{Reconstruction errors for the discontinuous function  $u_d$ on a
  uniform grid. The discontinuity is on a grid interface from N=640 onwards.}
\label{tab:eps:jumpnew}
\end{table}

\begin{table}
\begin{minipage}[t]{0.49\textwidth}
Smooth test, quasi-uniform grid\\
\footnotesize
\begin{tabular}{r|cc|cc|}
N & $\|E\|_1$ & rate & $\|E\|_\infty$ & rate \\
\hline
  20 &  5.63e-03 &      & 2.17e-02 &     \\ 
  40 &  8.65e-04 & 2.70 & 3.85e-03 & 2.50\\ 
  80 &  1.08e-04 & 3.01 & 5.26e-04 & 2.87\\ 
 160 &  1.26e-05 & 3.10 & 5.36e-05 & 3.30\\ 
 320 &  1.47e-06 & 3.09 & 7.00e-06 & 2.94\\ 
 640 &  1.75e-07 & 3.08 & 8.85e-07 & 2.98\\ 
1280 &  2.11e-08 & 3.05 & 1.11e-07 & 3.00\\ 
2560 &  2.59e-09 & 3.03 & 1.39e-08 & 3.00\\ 
\end{tabular}
\end{minipage}
\begin{minipage}[t]{0.49\textwidth}
Smooth test, random grid\\
\footnotesize
\begin{tabular}{r|cc|cc|}
N & $\|E\|_1$ & rate & $\|E\|_\infty$ & rate \\
\hline
  20 & 5.36e-03 &      & 1.87e-02 &     \\ 
  40 & 4.97e-04 & 3.43 & 1.35e-03 & 3.79\\ 
  80 & 5.07e-05 & 3.29 & 1.35e-04 & 3.32\\ 
 160 & 5.46e-06 & 3.22 & 1.67e-05 & 3.02\\ 
 320 & 6.25e-07 & 3.13 & 2.08e-06 & 3.01\\ 
 640 & 7.46e-08 & 3.07 & 2.64e-07 & 2.98\\ 
1280 & 9.11e-09 & 3.03 & 3.15e-08 & 3.07\\ 
2560 & 1.12e-09 & 3.02 & 4.19e-09 & 2.91\\ 
\end{tabular}

\end{minipage}
\caption{Reconstruction errors for the smooth function  $u_s$  on
  non-uniform grids. The quasi uniform grids have cell centers at 
  $x_j=j/N+0.1\sin(20\pi j/N)/10$ and the random ones at
  $x_j=j/N+0.25/N r_j$ where the $r_j$ are uniformly distributed in
  $[-0.5,0.5]$. }  
\label{tab:eps:nonunif}
\end{table}

\paragraph{One-dimensional reconstruction tests}
We set up $U_j$ with the cell averages of the smooth function $u_s(x)=\sin(2\pi x-\sin(2\pi x)/\pi)$ on the domain $[0,1]$ and of the discontinuous function $u_d(x)=e^{-x^2}+0.1H(x)$
on $[-0.5,0.5]$, where $H(x)$ is the Heavyside function. Then boundary extrapolated data $U_{j+1/2}^{\pm}$ are computed with the reconstruction and compared with limits of the function $u_s$ (resp. $u_d$) for  $x\to x_{j+1/2}^{\pm}$.

Table
\ref{tab:eps:smooth} compares the choice $\epsilon=h$ and the
classical choice $\epsilon=10^{-6}$ in the case of the smooth
function. It is clear that our choice yields lower errors at all grid
resolutions and a much more regular convergence pattern. The results
for $\epsilon=h^2$ (not reported) are in between the two.

%Table \ref{tab:eps:jump} is about the same test for the discontinuous
%function. In order to interpret the results, it is important to notice
%that the discontinuity is located exactly at a cell interface.  In
%this case, the choice $\epsilon=h^2$ (not shown) yields rates closer
%to those obtained with $\epsilon=10^{-6}$, but we remark that all
%convergence rates would obviously be degraded to $1$ and $0$ (resp. in the
%$1$- and the $\infty$-norm) if the discontinuity were not exactly on a
%cell interface.

Table \ref{tab:eps:jumpnew} is about the same test for the discontinuous function. Here we can observe two very different regimes, depending on whether the discontinuity is located exactly at a grid interface or not. In our tests the discontinuity is located  at $x=1/320$, so that it is located exactly at a cell interface only from N=640 onwards, while being inside a cell for coarser grids. The convergence rates are obviously degraded to $1$ and $0$ (resp. in the $1$- and the $\infty$-norm) when the discontinuity is not located exactly at a cell interface. We notice that both choices for $\epsilon$ yield very similar errors in both norms, until $N=320$. From $N=640$ onwards, since the discontinuity is located exactly at an interface, the errors improve greatly. In this regime, $\epsilon=10^{-6}$ can exploit more the favourable situation and produce very small errors. This situation is however quite rare in an evolutionary problem and thus we can conclude that the choice $\epsilon=h$ is a better one, on the grounds that it is more accurate for the smooth parts and comparable to $\epsilon=10^{-6}$ around discontinuities.

Finally, in Table \ref{tab:eps:nonunif}, the CWENO reconstruction with
$\epsilon=h_j$, i.e.\ the local mesh size, is tested on non-uniform grids, both of quasi-regular and
random type (see the caption of Table \ref{tab:eps:nonunif} for the definition of these grids), showing remarkable robustness in the order of convergence
and errors very close to those obtained on uniform grids, suggesting that the choice $\epsilon\propto h$ is a good one, and that the reconstruction algorithm works very well even on non-uniform grids.

\paragraph{Two-dimensional reconstruction tests}
We apply the reconstruction to the exact cell averages
of the smooth function
$u^\text{exa}(x,y)=\sin(\pi x)\cos(\pi y)$ on the unit square, with periodic boundary conditions.
The error $E$ is computed by comparing the reconstructed function $u^\text{rec}$ with
 $u^\text{exa}$ on a uniform grid of reference $\mathcal{G}_\text{ref}$, in both the $1$-norm and the
$\infty$-norm. The reconstruction $u^\text{rec}$ is defined by $u^\text{rec}|_{\Omega_j} = P$, for $j=1,\ldots,N$, with $P$ being the reconstruction polynomial \eqref{finalrec1d} and \eqref{finalrec} on the cell $\Omega_j$.
Cell averages $U_j$ are set using the exact computation of the integral, namely $U_j = \int_{\Omega_j} u^\text{exa} d\Omega$.
 In details:
\begin{equation*} %\label{errorFormula}
 \left\| E \right\|_{q} = \left( \frac{ \sum_{G \in \mathcal{G}_{ref}} \left| u^\text{rec}(G)-u^\text{exa}(G) \right|^q}{|\mathcal{G}_\text{ref}|} \right)^{1/q}, \; \; \; q=1,\infty.
\end{equation*}
The reference grid $\mathcal{G}_\text{ref}$ is a uniform grid such that each cell is finer than the smallest cell allowed by the adaptive algorithm. Although it is sufficient to have an accurate reconstruction only on the boundary of each cell (since the reconstruction is needed to compute numerical fluxes on quadrature points, see Fig.\ \ref{fig:gausspoints}), by choosing this grid we show that the method indeed provides a uniform accuracy all over the domain.

 Table
\ref{tab:errunif:smooth} reports the reconstruction errors observed on
uniform grids. Next, an h-adapted grid with $3$ levels of refinement was generated
by recursively refining the cells of a $N_\text{START}\times N_\text{START}$ uniform grid
where $\|P^1_{\text{OPT}}-P_{\text{OPT}}\|_2^2>0.01 h_j^2$, 
with $ P^1_{\text{OPT}}$ being the optimal first degree 
polynomial computed with the same least square procedure of $P_{\text{OPT}}$.
We choose $N_\text{START}=8$ and denote this grid as $\mathcal{G}_0$. From this grid we generate grids $\mathcal{G}_k$, for $k=1,\ldots$, by subdividing each cell of $\mathcal{G}_0$ into $4^k$ equal cells.
For each grid $\mathcal{G}_k$, $k=0,\ldots$, we define $N_\text{C}=N_\text{START} \cdot 2^k = 2^{3+k}$.
The result of the convergence test on
non-uniform meshes are reported in Table
\ref{tab:erradap:smooth}. The two-dimensional results are in line with those obtained in one dimension, confirming the good quality of the choice $\epsilon\propto h$.

%\matteo{Armando, come sono calcolati gli errori? Per la norma-$\infty$
%  immagino si tratti dei valori ai bordi: dove esattamente? Nei nodi
%  di quadratura? E la norma-1 a cosa si riferisce?}
%\matteo{Armando, puoi fare anche un test con una soluzione
%  discontinua? Per esempio $\sin(2\pi x-\sin(2\pi x)/\pi)+\chi(x,y)$
%  dove $\chi$ \`e la funzione caratteristica di un cerchio di centro
%  $(0.5,0.5)$ e raggio circa $1/4$. Il raggio andrebbe scelto in modo
%  che rimangano delle celle medie dentro al cerchio. Dovrai usare
%  pi\`u di 3 livelli, immagino}

\begin{table}
\begin{minipage}[t]{0.49\textwidth}
Uniform grid, $\epsilon=10^{-6}$\\
\footnotesize
\begin{tabular}{r|cc|cc|}
N & $\|E\|_1$ & rate & $\|E\|_\infty$ & rate \\
\hline
$  8^2$ &   7.19e-02  &       & 3.11e-01 & \\
$  16^2$ &   1.36e-02  &  2.41 & 9.61e-02  &  1.69 \\
$  32^2$ &   2.16e-03  &  2.65 & 2.51e-02  &  1.94 \\
$ 64^2 $&   2.72e-04  &  2.99 & 4.11e-03  &  2.61 \\
$ 128^2$ &   2.68e-05  &  3.34 & 4.33e-04  &  3.25
\end{tabular}
\end{minipage}
\hfill
\begin{minipage}[t]{0.49\textwidth}
Uniform grid, $\epsilon=h_j$\\
\footnotesize
\begin{tabular}{r|cc|cc|}
N & $\|E\|_1$ & rate & $\|E\|_\infty$ & rate \\
\hline
$  8^2$ & 5.32e-02  &        & 1.90e-01  &  \\
$  16^2$ & 8.51e-03  &  2.64  & 3.74e-02  &  2.34\\
$  32^2$ & 9.56e-04  &  3.15  & 3.79e-03  &  3.30\\
$ 64^2$ & 9.45e-05  &  3.34  & 2.55e-04  &  3.89\\
$ 128^2$ & 1.09e-05  &  3.12  & 3.20e-05  &  3.00
\end{tabular}
\end{minipage}
\caption{Reconstruction errors for a smooth function on a uniform
  grid. The average order of accuracy is $2.84$ for $\|E\|_1$ and $2.35$
  for $\|E\|_\infty$ with $\epsilon=10^{-6}$, and $3.10$ for $\|E\|_1$
  and $3.23$ for $\|E\|_\infty$ with $\epsilon=h$.}
\label{tab:errunif:smooth}
\end{table}

\begin{table}
\begin{minipage}[t]{0.49\textwidth}
Adaptive grid, $\epsilon=10^{-6}$\\\
\footnotesize
\begin{tabular}{r|cc|cc|}
$N_{\text{C}}$ & $\|E\|_1$ & rate & $\|E\|_\infty$ & rate \\
\hline
$  8^2$ &   6.19e-02  &       & 3.11e-01 & \\
$  16^2$ &   1.11e-02  &  2.48 & 9.61e-02  &  1.69 \\
$  32^2$ &   1.75e-03  &  2.67 & 2.53e-02  &  1.93 \\
$ 64^2$ &   2.33e-04  &  2.91 & 6.40e-03  &  1.98
\end{tabular}
\end{minipage}
\hfill
\begin{minipage}[t]{0.49\textwidth}
Adaptive grid, $\epsilon=h_j$\\
\footnotesize
\begin{tabular}{r|cc|cc|}
$N_{\text{C}}$ & $\|E\|_1$ & rate & $\|E\|_\infty$ & rate \\
\hline
$  8^2$ &   4.48e-02  &       & 1.90e-01 & \\
$  16^2$ &   7.00e-03  &  2.68 & 3.74e-02  &  2.34 \\
$  32^2$ &   8.16e-04  &  3.10 & 3.79e-03  &  3.30 \\
$ 64^2 $&   8.91e-05  &  3.20 & 2.55e-04  &  3.89
\end{tabular}
\end{minipage}
\caption{Reconstruction errors for a smooth function on an adaptive
  grid. The first column shows the number of cells on the coarsest
  grid. For each test we choose $L=3$ levels of refinement. The average
  order of accuracy is $2.68$ for $\|E\|_1$ and $1.87$ for
  $\|E\|_\infty$ with $\epsilon=10^{-6}$, and $3.00$ for $\|E\|_1$ and
  $3.19$ for $\|E\|_\infty$ with $\epsilon=h$.}
\label{tab:erradap:smooth}
\end{table}

\paragraph{Convergence tests on adaptive grids: notation}
The following tests will compare the numerical schemes on uniform and adaptive grids. We will consider a second order scheme (minmod reconstructions, Heun timestepping, numerical entropy evaluated on the cell averages) and the third order scheme described in this paper (CWENO reconstruction, three-stage third order SSP-RK and numerical entropy indicator  computed with Gauss quadratures with two points per direction). The order will be denoted in the legends with the letter $r$.

In all tests the numerical fluxes are the Local Lax Friedrichs ones, i.e.
\[
F(u,v) = \frac12\left[f(u)+f(v)-\alpha(u-v)\right]
\]
where $\alpha$ is the largest eigenvalue among those of $f'(u)$ and $f'(v)$. The numerical entropy fluxes are chosen accordingly  as
\[
\Psi(u,v) = \frac12\left[\psi(u)+\psi(v)-\alpha(\eta(u)-\eta(v)\right]
\]
with the same value for $\alpha$ (see \cite{PS:2011:entropy}). For the computation of the numerical entropy production we have employed $\eta(u)=u^2$ for scalar tests and the physical entropy for the gas-dynamics tests.

The convergence tests are conducted by performing sequences of computations with a coarse grid of size $N_0=M2^k$ for $k=0,\ldots,K$ for some integer $M$.
We say that a simulation employs $\ell$ cell levels if the scheme is allowed to chose among $\ell$ different cell sizes, i.e. $H,H/2,\ldots,H/2^{\ell-1}$, where $H$ is the size of the uniform coarse mesh. Let $\ell(k)$ be the number of levels employed in the $k$-th test of a sequence. In the legends of the convergence graphs, the number of levels is indicated by an integer $\mathtt{L}$ if $\ell(k)=\mathtt{L}$, by $\mathtt{L+}$ (respectively $\mathtt{L++}$) if $\ell(k)=\mathtt{L}+k$ (respectively $\ell(k)=\mathtt{L}+2k$). Obviously $\mathtt{L}=1$ corresponds to uniform grids. The refinement threshold is set to $S_{\text{ref}} =S_0\mathrm{s}^{-k}$ in the $k$-th computation for some scaling factor $s\geq1$, where $S_0$ is the threshold used in the coarsest computation of the sequence. 
The choice of the parameter $s$ should be guided by the features of the solution on which one wants h-adaptivity to act.
Since the numerical entropy production scales as $O(h^r)$ on smooth flows (see Section 3), we condict the tests using $s=2^r$ if the solution is smooth everywhere and h-adaptivity is expected to distinguish between high and low frequency regions. On the other hand, the numerical entropy production scales as $O(1/h)$ on shocks, $O(1)$ on contacts and $O(h)$ on corner points and thus we take $s=2$ on flows with singularities so that h-adaptivity acts on rarefaction corners and stronger discontinuities. These choices for $s$ ensure that the level of refinement induced by the choice of $S_0$ on the first grid is replicated in the $k$-th test.
Since in h-adaptive tests the number of cell is variable during the evolution, in the graphs for the convergence tests we employ the time-average of the number of cells employed in each timestep.

\paragraph{One-dimensional linear transport test (1)}
We evolve the initial data $u(x,0)=\sin(\pi x-\sin(\pi x)/\pi)$ on $[-1,1]$
up to $t=1.0$ with the linear transport equation $u_t+u_x=0$ and
periodic boundary conditions.

\begin{figure}
\begin{center}
\includegraphics[width=0.7\textwidth]{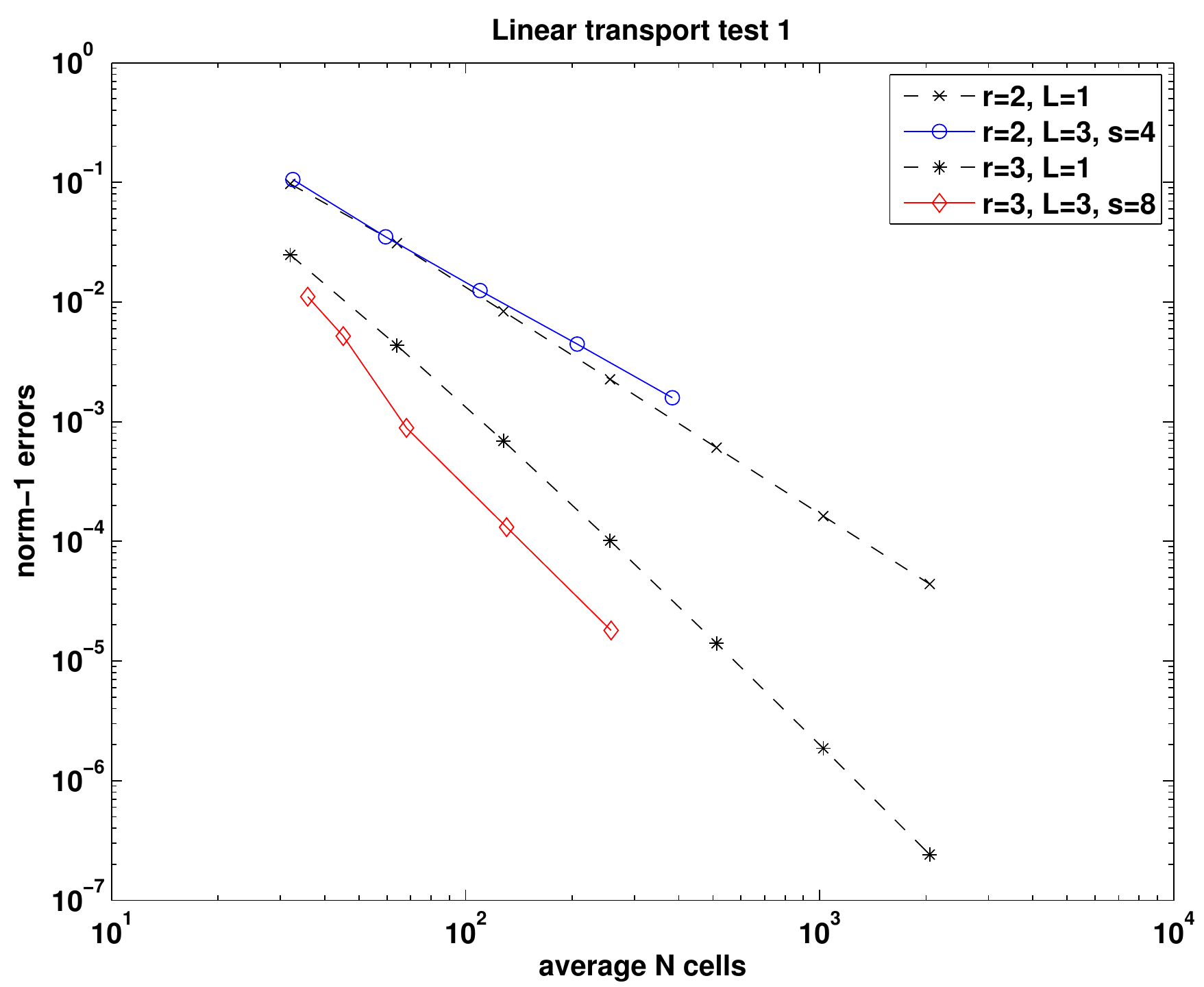}
\end{center}
\caption{Norm-1 errors for uniform and adaptive grids for the linear transport test (1). The number of levels is kept fixed, while the coarse cell size and  the refinement threshold are progressively lowered. The symbols in the legend refer to the order of the scheme ($r$), the number of different cell sizes ($\mathtt{L}$) and the scaling of the refinement threshold between runs ($s$).}
\label{fig:lintra1}
\end{figure}

The purpose of this test is to ensure that, on a single scale smooth solution, with a proper choice of the parameter $s$, the behavior of the adaptive scheme is similar to that of a uniform grid method.

In the set of runs, the number of grid refinement levels was kept at $\mathtt{L}=3$, while the coarse grid size and the refinement threshold have been diminished (see Fig.\ \ref{fig:lintra1}). The refinement threshold for the coarsest run in the series ($S_0$) is taken to be $10^{-1}$ for the second order scheme and $10^{-2}$ for the third order scheme. As it appears from the figure, the behavior of the adaptive scheme with the choice $s=2^r$ is similar to the one given by the uniform grid method. We observe a slight improvement of the adaptive scheme for $r=3$, while for $r=2$ there is really no advantage in using an adaptive mesh in such a simple test case.

\begin{figure}
\begin{center}
\includegraphics[width=0.7\textwidth]{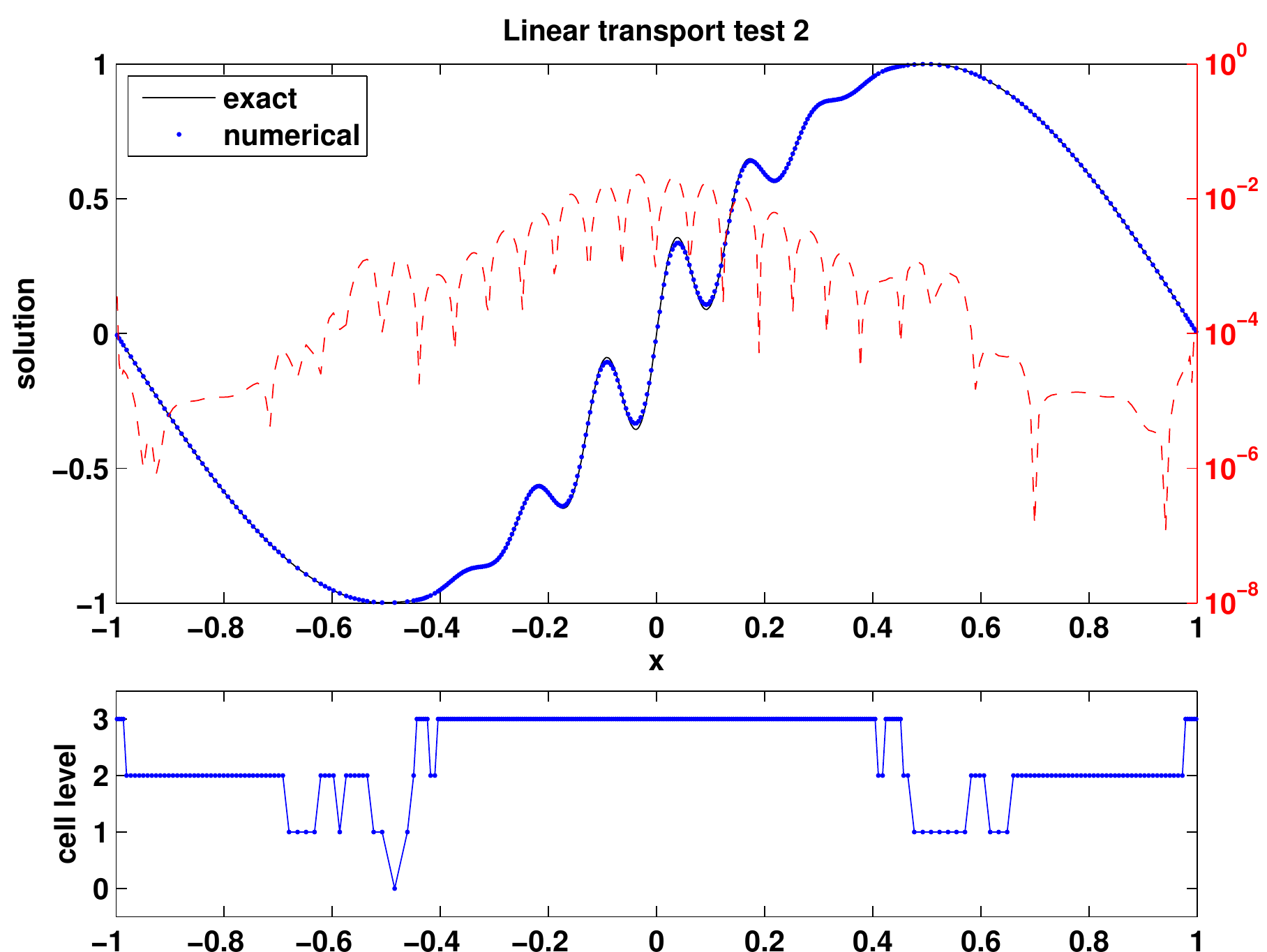}
\\
\includegraphics[width=0.7\textwidth]{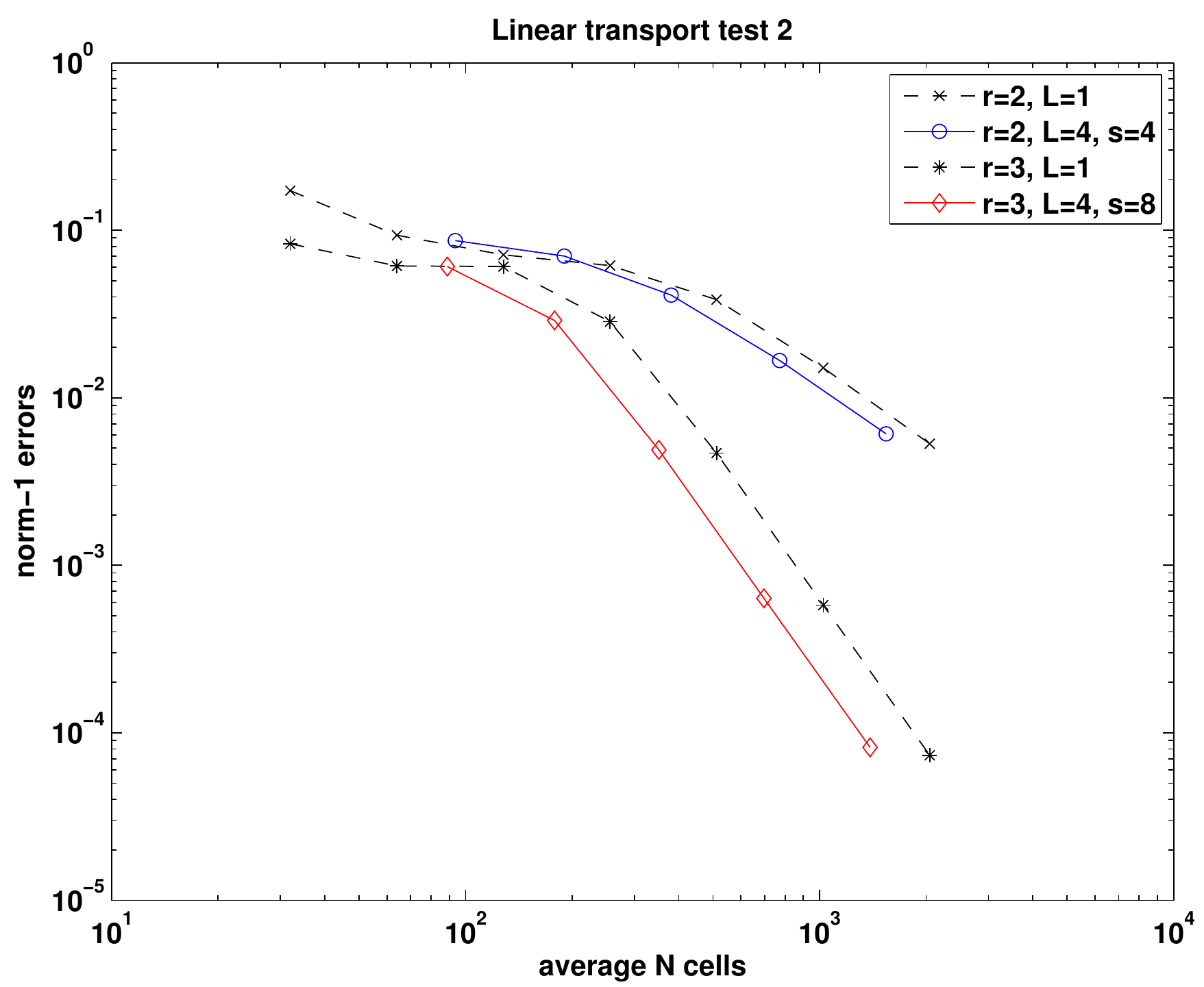}
\end{center}
\caption{Norm-1 errors for uniform and adaptive grids for the linear transport test (1). On the left, an example of a solution (top) with the grid levels employed (bottom). In the top panel the dashed red line is the error of the numerical solution and should be read against the logarithmic scale on the right vertical axis. On the right, the convergence test.}
\label{fig:lintra2}
\end{figure}

\paragraph{One-dimensional linear transport test (2)}
Next, we want to test the situation in which the adaptive scheme may compute on a finer grid the rapidly varying regions of a smooth solution. To this end, we evolve the initial data $u(x,0)=\sin(\pi x)+\tfrac14\sin(15\pi x)e^{-20x^2}$ on $[-1,1]$
up to $t=2.0$ with the linear transport equation $u_t+u_x=0$ and
periodic boundary conditions (see Fig.\ \ref{fig:lintra2}, left). The top graph shows an example of  numerical solution obtained with the adaptive scheme ($N_0=64$,$\mathtt{L}=3$,$S_{\text{ref}}=2\times10^{-5}$), while the bottom one shows the cell levels in use at final time, revealing that small cells are selected in the central region where the solution contains higher frequencies. The dashed red line in the top graph is the error of the adaptive solution and should be read against the logarithmic scale on the right. It reveals that no spurious oscillations or other defects occur at grid discontinuities.

The results of the convergence tests are shown on the right in Fig.\ \ref{fig:lintra2}. First, note that the runs on uniform grids cannot resolve the high frequencies in the central area and thus they do not show convergence until $N_0\geq 128$. For example, the convergence rates estimated from the errors shown in the figure for the third order scheme are 0.44, 0.01, 1.09, 2.61, 3.02 and 2.98.

For the adaptive grid tests we used $N_0=16\cdot2^k$ and $\mathtt{L}=4$, so that even in the coarsest run the scheme was able to use some cells of size $1/128$-th of the domain. $S_0$ was set to $0.005,0.001$ for the schemes of order $\mathrm{r}=2,3$ respectively. These values were chosen in such a way that, for $k=0$, only cells in the region of the 4 central peaks ($[-0.2,0.2]$ in the initial/final data shown in the left panel of the figure) had a chance to be refined to the 4-th level. With respect to the previous test, obviously here an adaptive scheme has much room for improvement over uniform grid ones and all tested sequences show improvements over the uniform case.

Since the solution is smooth everywhere, we expect that scaling the refinement threshold with $s=8$ in the scheme with $r=3$ would cause every run to use the same refinement pattern. As a matter of fact this choice leads to the best convergence rates (1.06, 2.62, 2.99, 2.97). 

\paragraph{A two-dimensional scalar test with smooth solution}
Next we consider a scalar problem in two space dimensions with a
smooth solution which is known in closed form, originally presented in
\cite{Davies:85}. It is sometimes
referred to as ``atmospheric instability problem'', as it mimics the
mixing of a cold and hot front. The equation is
\begin{equation}  \label{eq:swirl}
u_t + \nabla\cdot(v(x,y) u) =0 
\quad 
\text{ on } \Omega=[-4,4]^2
\end{equation}
with initial data
$ u_0(x,y) = -\tanh(y/2)$
and velocity
\[ v(x,y) 
= \left[ -\tfrac{y}{r}\tfrac{f}{0.385},
\tfrac{x}{r}\tfrac{f}{0.385}\right] ,
\]
where $r=\sqrt{x^2+y^2}$ and $f=\tanh(r)/(\cosh(r))^2$. The exact
solution at $t=4$ is
\[ 
u(4,x,y)=-\tanh\left(
\tfrac{y}{2}\cos\left(\tfrac{4f}{0.385r}\right)
-\tfrac{x}{2}\sin\left(\tfrac{4f}{0.385r}\right)
\right)
\]

\begin{figure}
\begin{center}
\includegraphics[width=0.49\textwidth]{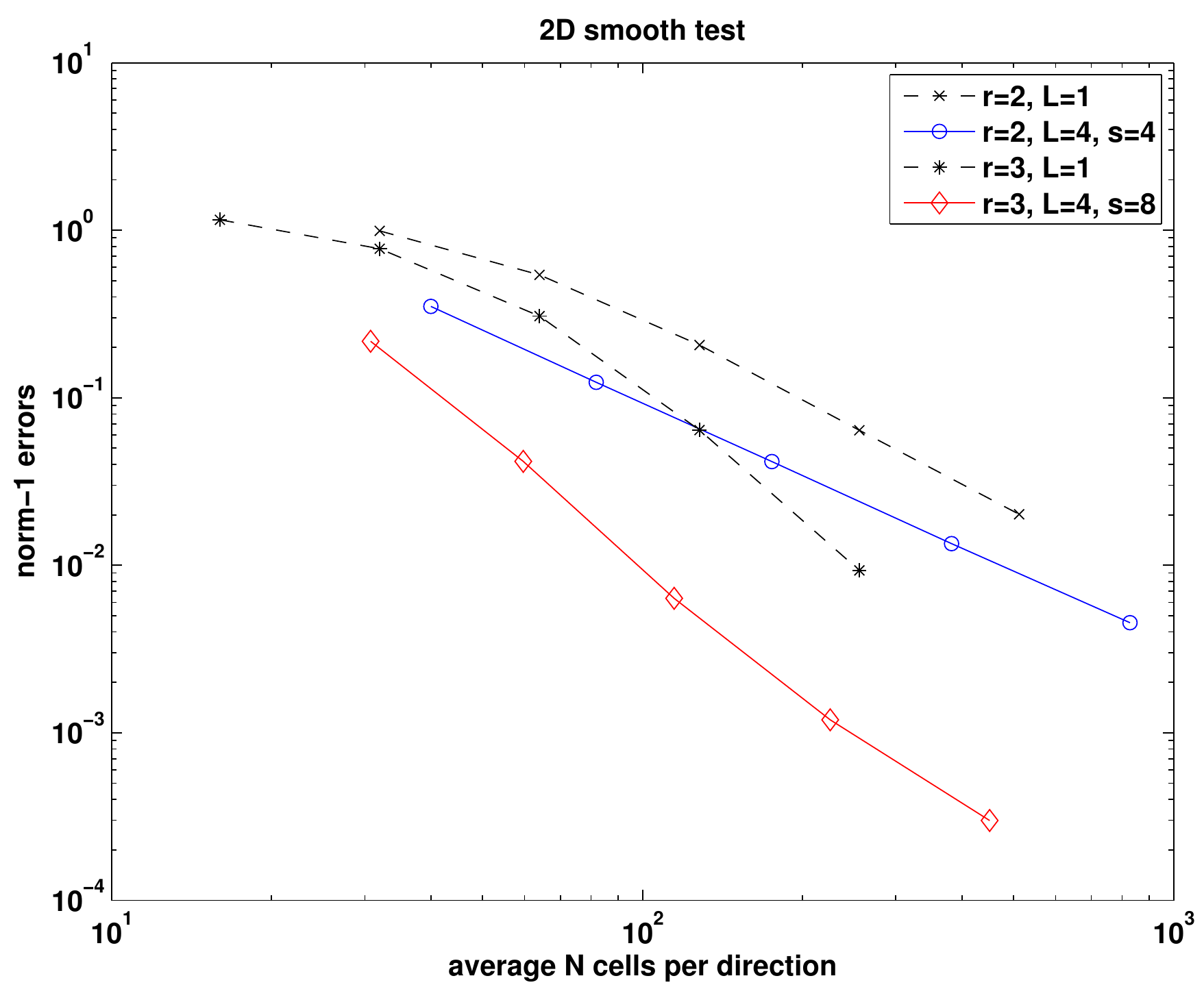}
\includegraphics[width=0.49\textwidth]{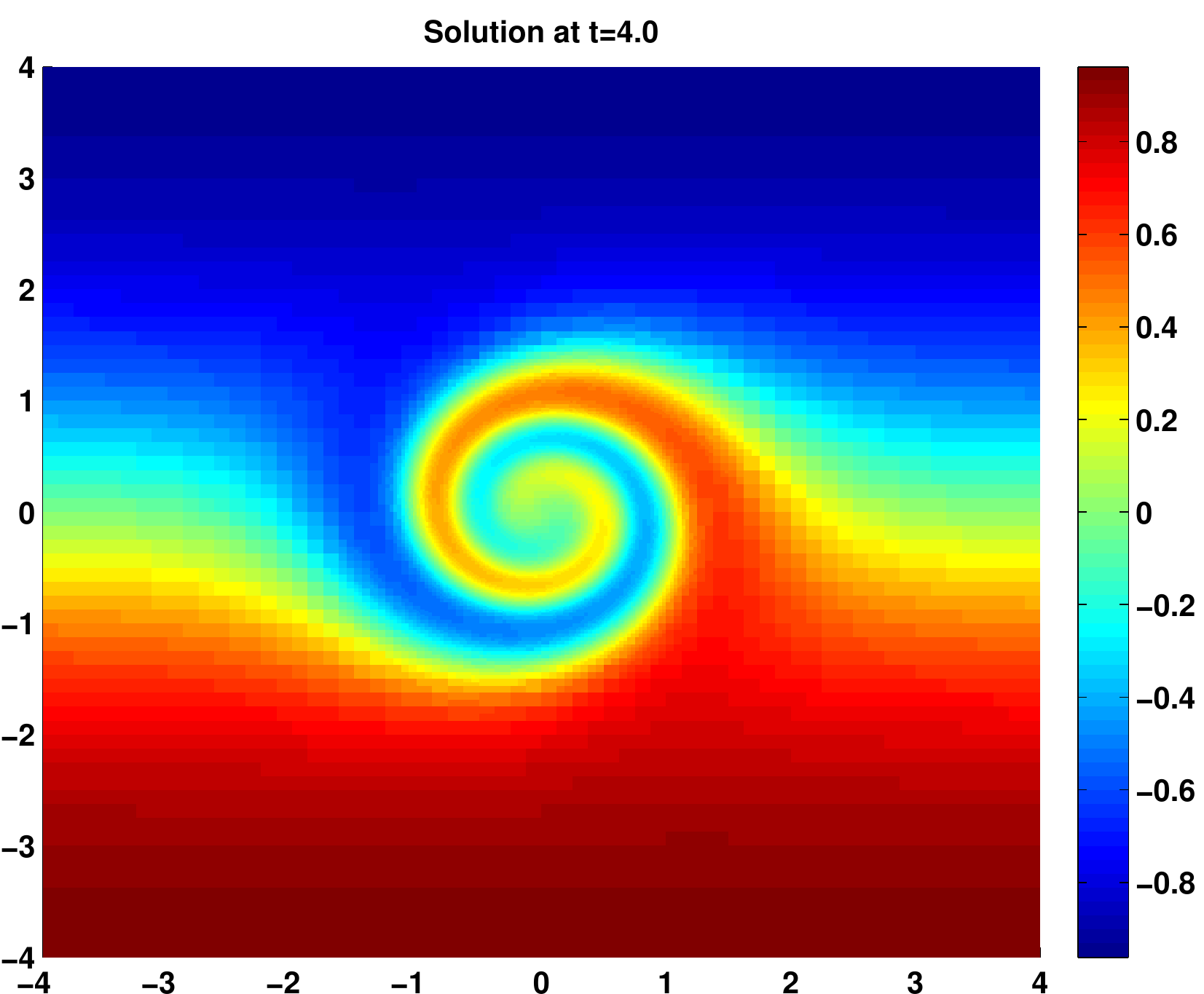}

\includegraphics[width=0.49\textwidth]{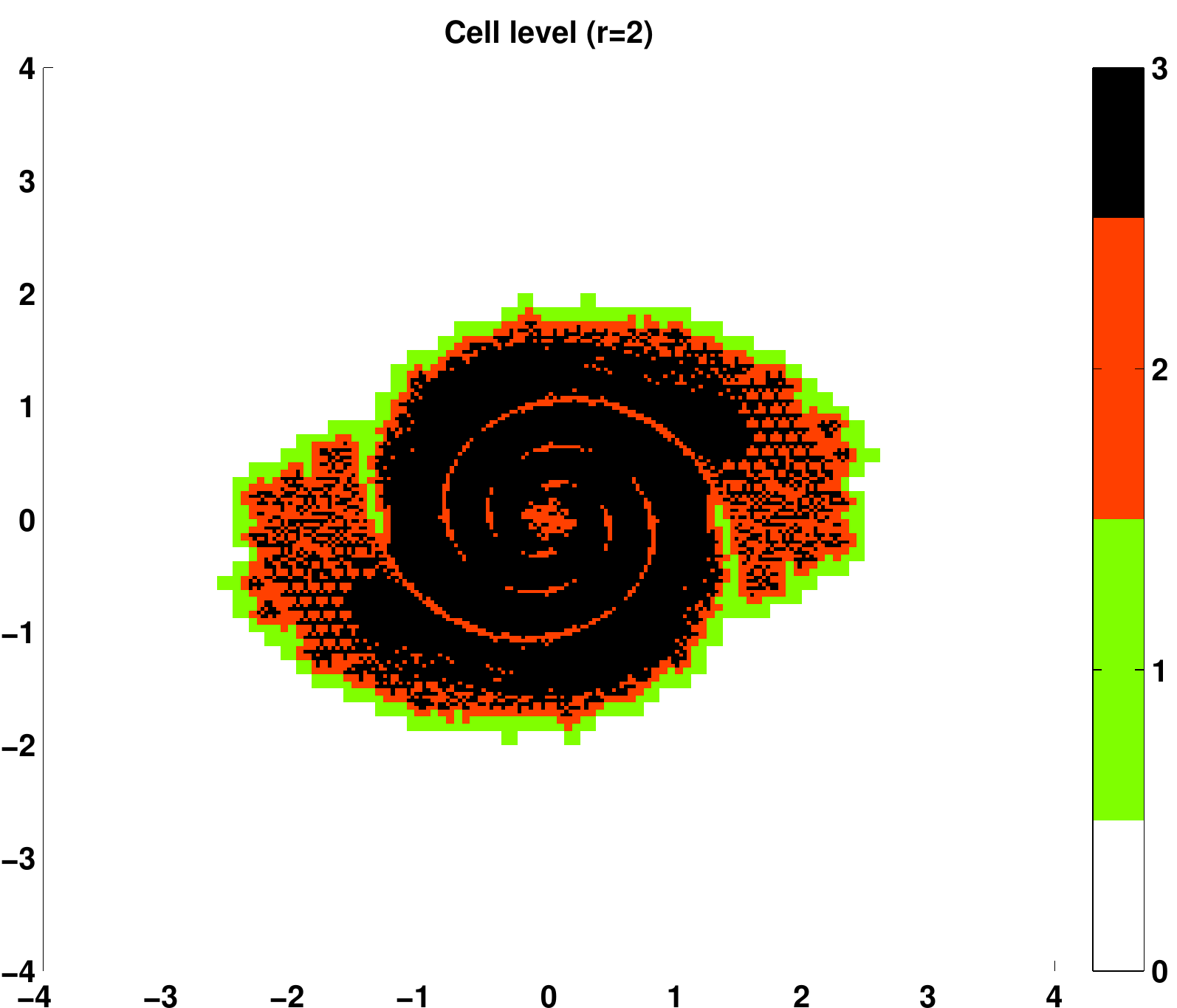}
\includegraphics[width=0.49\textwidth]{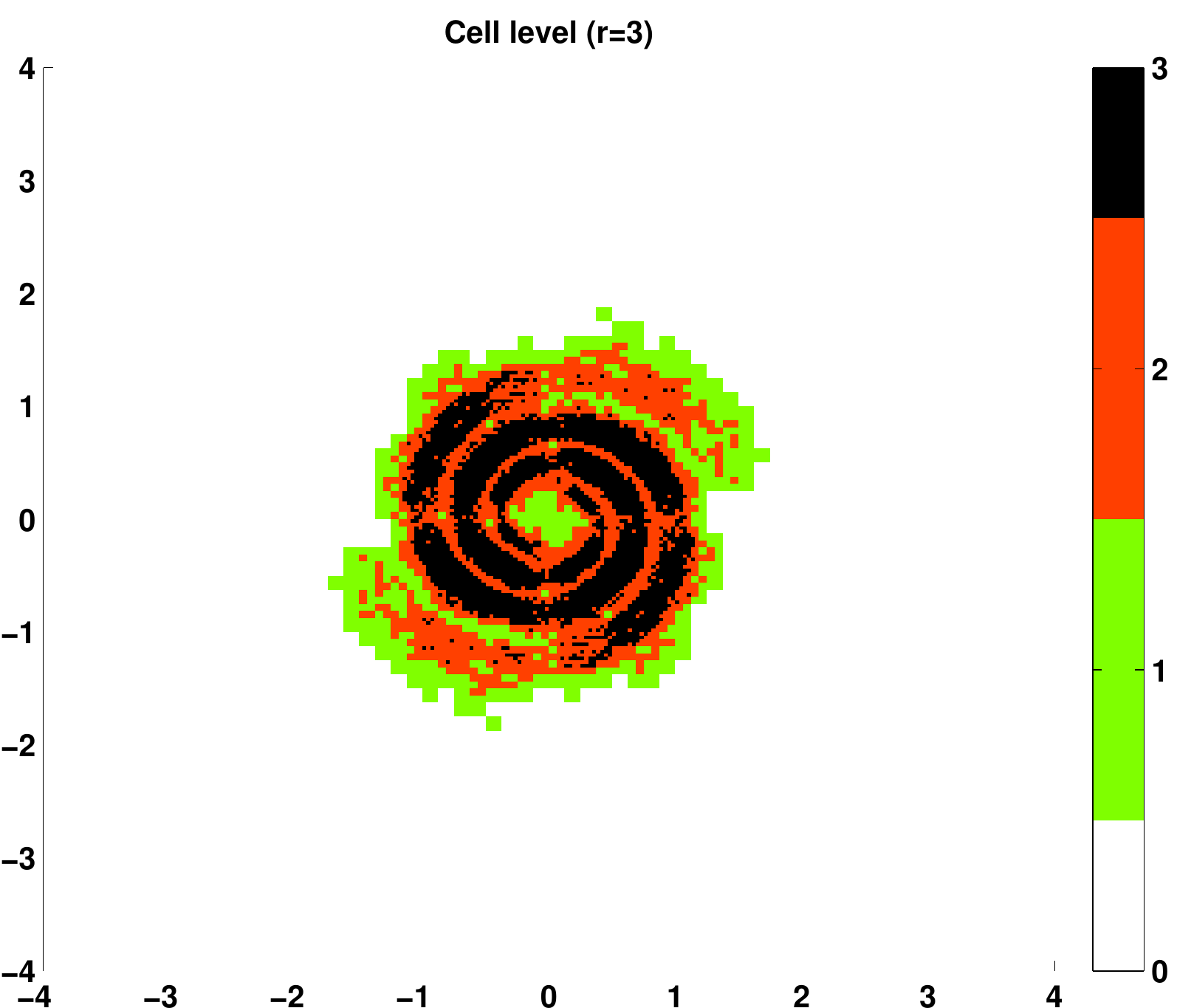}
\end{center}
\caption{Results of the 2d smooth test of equation \eqref{eq:swirl}. Norm-1 errors (top-left), solution at final time (top-right), grid refinement level at final time with the second and third order scheme with $N_0=64$, $\mathtt{L}=4$, $S_{\text{ref}}=6 \times 10^{-4}$ (bottom).}
\label{fig:swirl}
\end{figure}

Fig.\ \ref{fig:swirl} shows the solution at final time and the results of the convergence test conducted with $M=16$, $S_0=10^{-2}$ for both $s=2,3$. It is clear that third order shemes achieve far better efficiency, by computing solutions with lower errors and using fewer cells. Note that even the uniform grid third order scheme outperforms the second order adaptive scheme, which yields lower errors but also slightly lower convergence rates than the second order uniform grid scheme. The third order adaptive scheme shows a lower error constant than the corresponding uniform grid one.

The most striking difference between second and third order schemes, however, is the different refinement pattern employed. The bottom graphs of Fig.\ \ref{fig:swirl} show the meshes at final time for both schemes in the run with $N_0=64$, $\mathtt{L}=4$, $S_{\text{ref}}=6 \times 10^{-4}$. Clearly the third order scheme is able to employ small cells in much narrower regions, thus leading to additional computational savings.

\begin{figure}
\begin{center}
\includegraphics[width=0.49\linewidth]{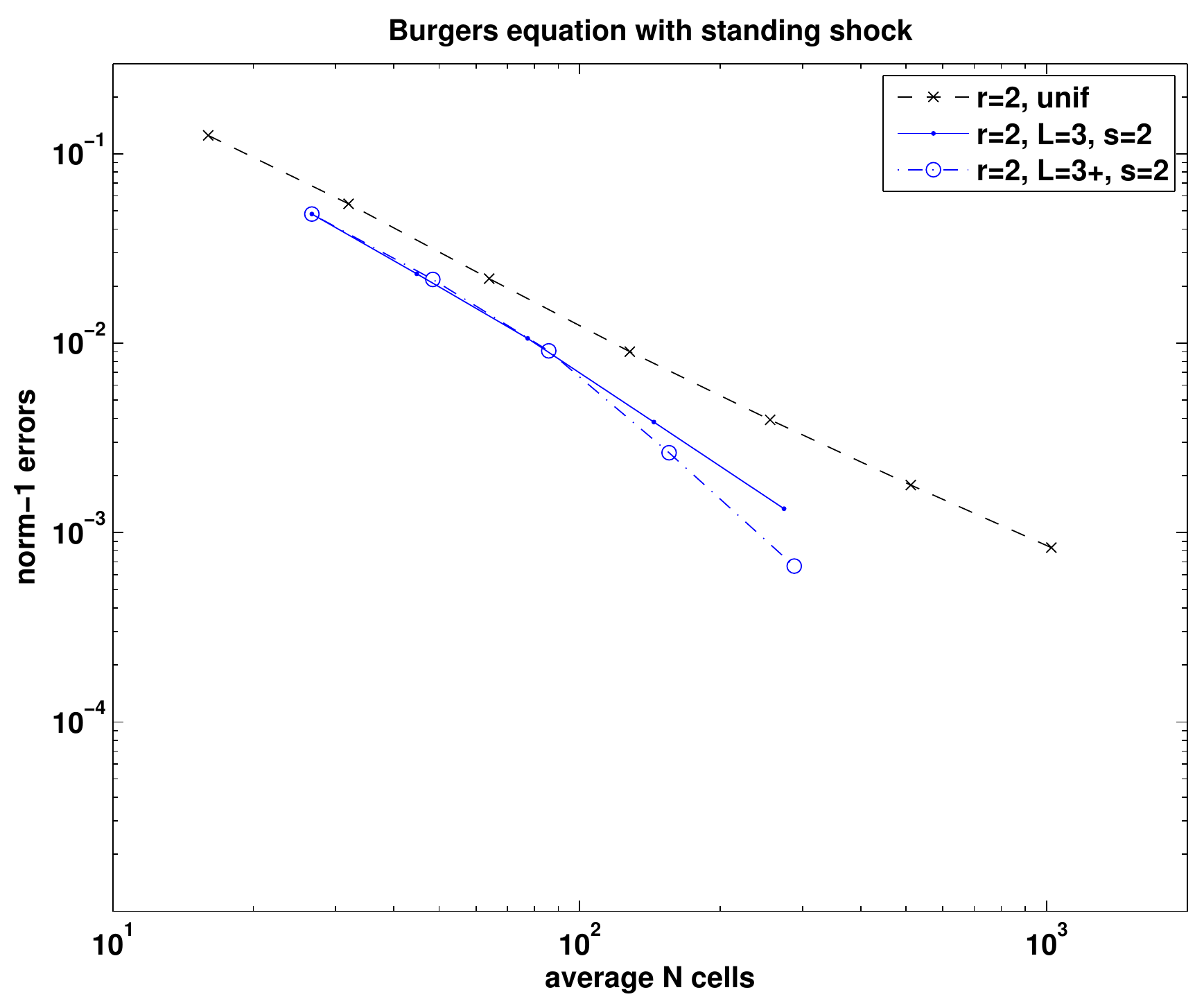}
\hfil
\includegraphics[width=0.49\linewidth]{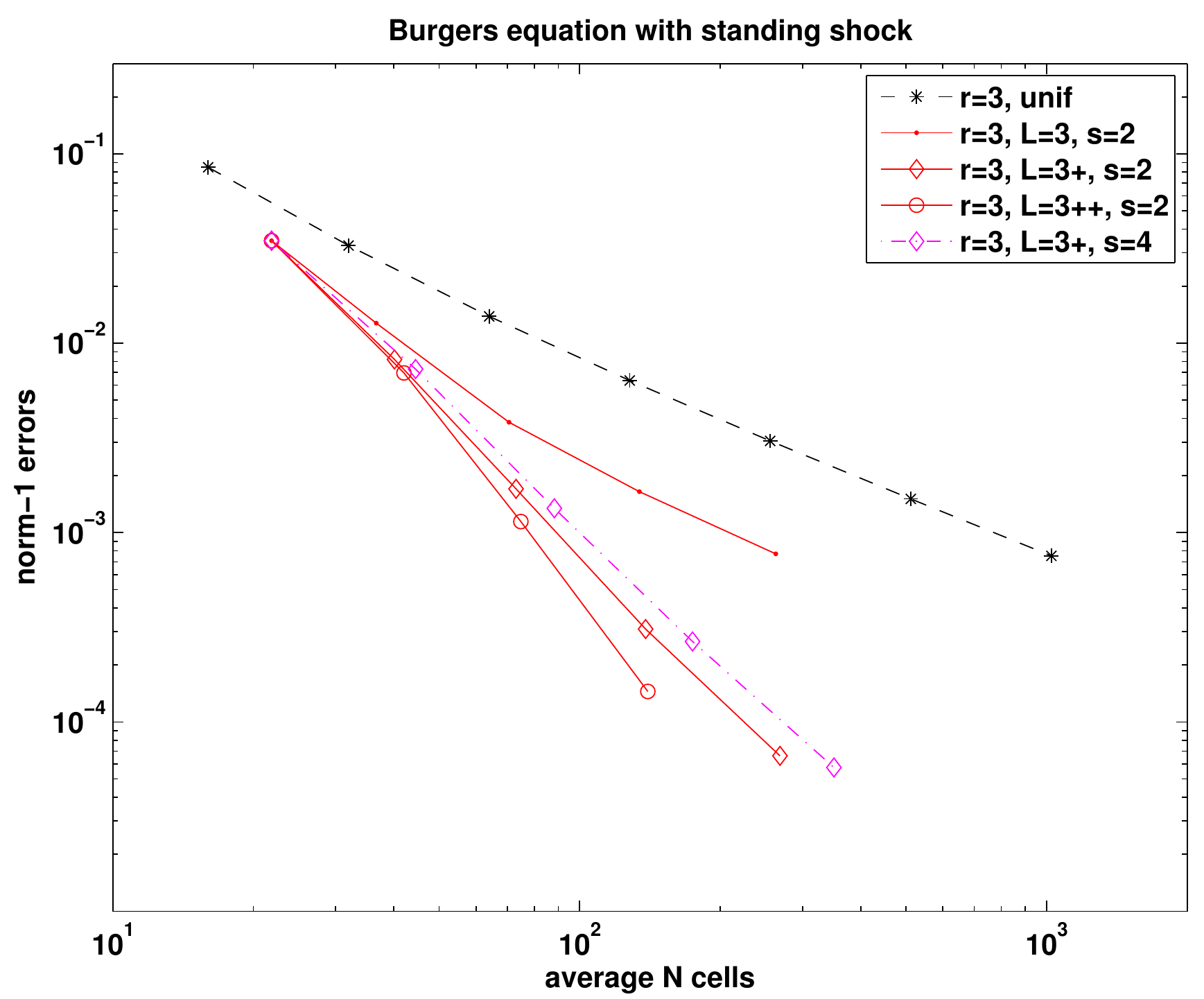}
\end{center}
\caption{Burgers equation after shock formation ($t=0.35$): error versus number of cells. Second order (left) and third order (right) scheme. The errors were computed by comparing with the exact solution computed with characteristics.}
\label{fig:burgersstanding}
\end{figure}

\begin{figure}
\begin{center}
\includegraphics[width=0.49\linewidth]{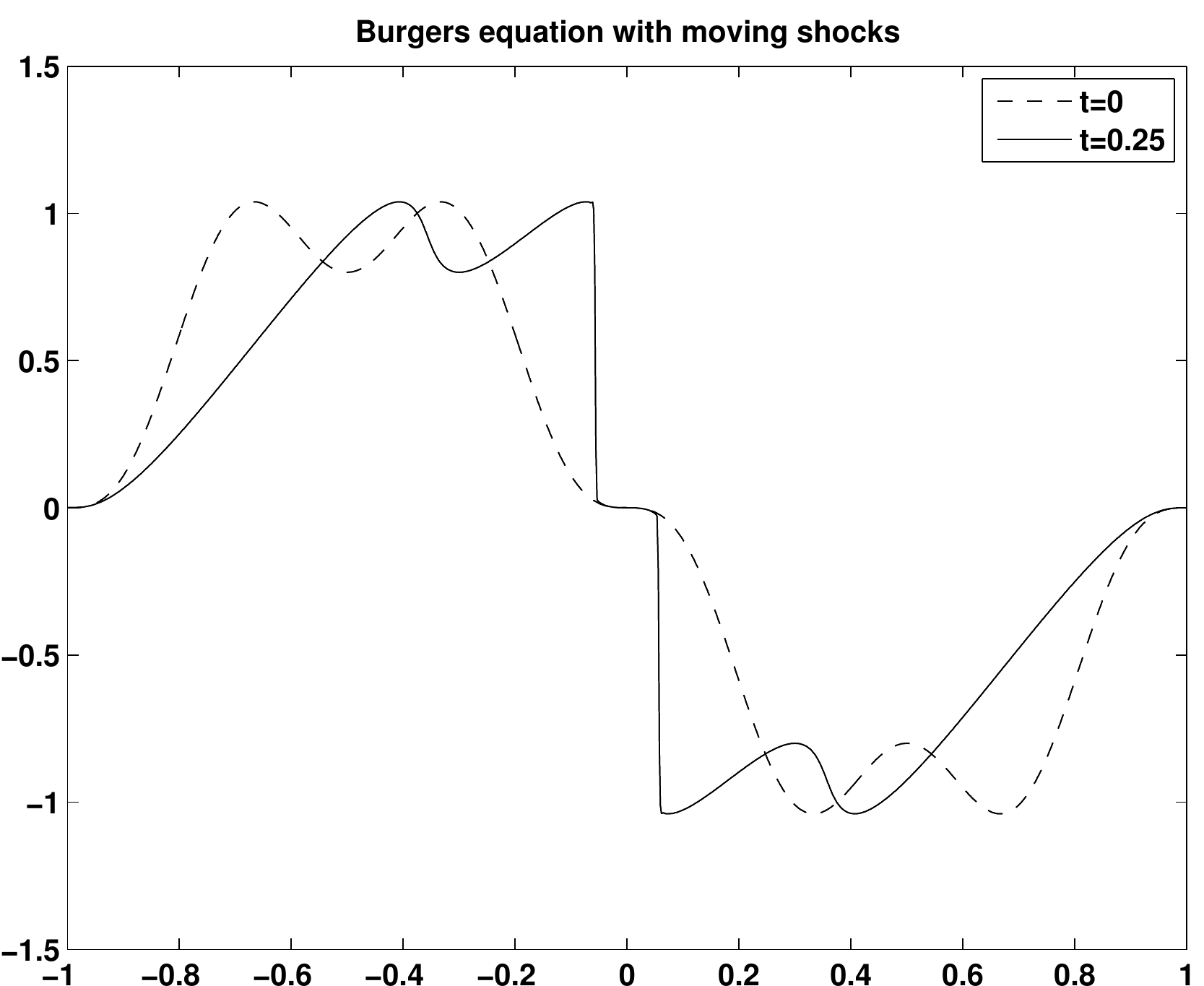}
\hfil
\includegraphics[width=0.49\linewidth]{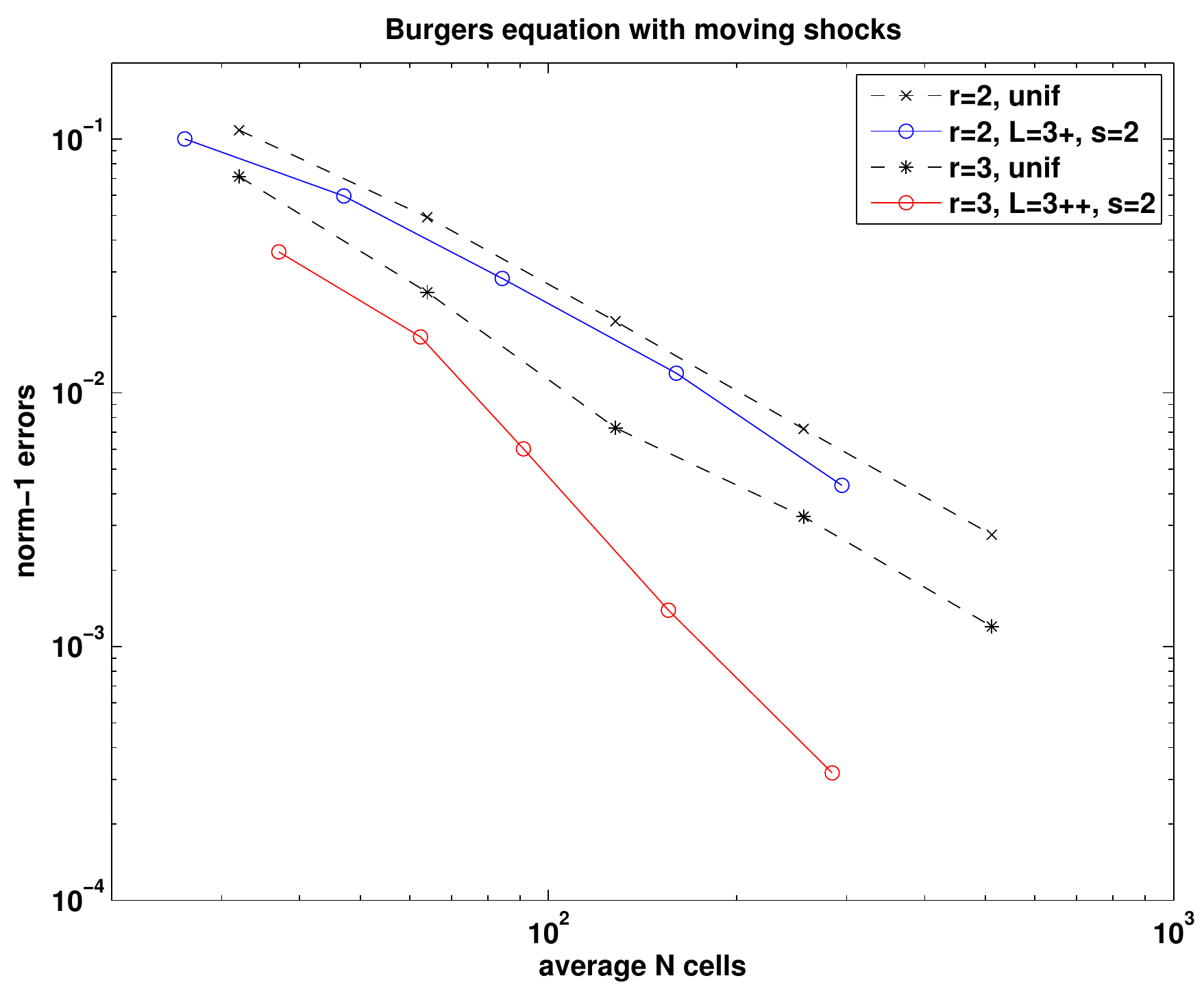}
\end{center}
\caption{Burgers equation with moving shocks: solution (left) and errors (right), obtained by comparing with a reference solution computed by the third order scheme on a uniform grid composed of $262144$ cells.}
\label{fig:burgersmoving}
\end{figure}

\paragraph{One dimensional Burgers equation} The first tests involving shocks employ Burgers' equation
\begin{equation}
\label{eq:burgers}
\partial_t u + \partial_x (\tfrac12 u^2) =0 
\end{equation}
on the domain $[-1,1]$ with periodic boundary conditions.
The results in Fig.\ \ref{fig:burgersstanding} are relative to the initial data $u_0(x)=-\sin(\pi x)$ and final time $t=0.35$, that gives rise to a standing shock located at $x=0$, while those
of Fig.\ \ref{fig:burgersmoving} refer to the initial condition 
$u_0(x)=-\sin(\pi x)+0.2\sin(5\pi x)$, that gives rise to two moving shocks.

In the error plot of Fig.\ \ref{fig:burgersstanding} first we note that the schemes with uniform meshes tend to become first order accurate, since most of the error is associated to the shock. For the adaptive runs we chose $M=16$ and $S_0=10^{-2}$ for both $r=2,3$. Adaptive schemes with a fixed number of levels do not improve much with respect to the uniform mesh ones: they give a smaller error, asymptotically approaching first order convergence, although with an error $2^{\mathtt{L}-1}$ times smaller than a uniform grid with the same number of cells. Allowing more levels when refining the coarse mesh yields an experimental order of convergence (EOC) of 2.0 (left). The third order scheme gives an EOC of 2.7 when increasing $\mathtt{L}$ by two every time the coarse mesh is refined (right). Note that the scaling factor for the refinement threshold can be taken as $s=2$, since here one is interested in refining only in the presence of the shock. For comparison purposes, we also show the result of using $s=4$ in the third order scheme, which gives a poorer work-precision performance due to excessive refinement.

Fig.\ \ref{fig:burgersmoving} refers to a more complex situation, since in this case the shocks are moving (and thus also coarsening comes into play) and there is a richer smooth structure away from the shocks. Here the advantage of third order schemes is more clear, even on uniform grids. This better performance is confirmed by the adaptive codes ($M=16$, $S_0=0.5$ for the second order case and $S_0=0.1$ for the third order one) and the situation is quite similar to the smooth test cases, except that the best convergence rates are obtained, similarly to the previous test with a shock, with scaling factor $s=2$ instead of $8$ (the EOC of the adaptive schemes are $1.7$ and $2.6$).

\begin{figure}
\begin{center}
\includegraphics[width=0.49\linewidth]{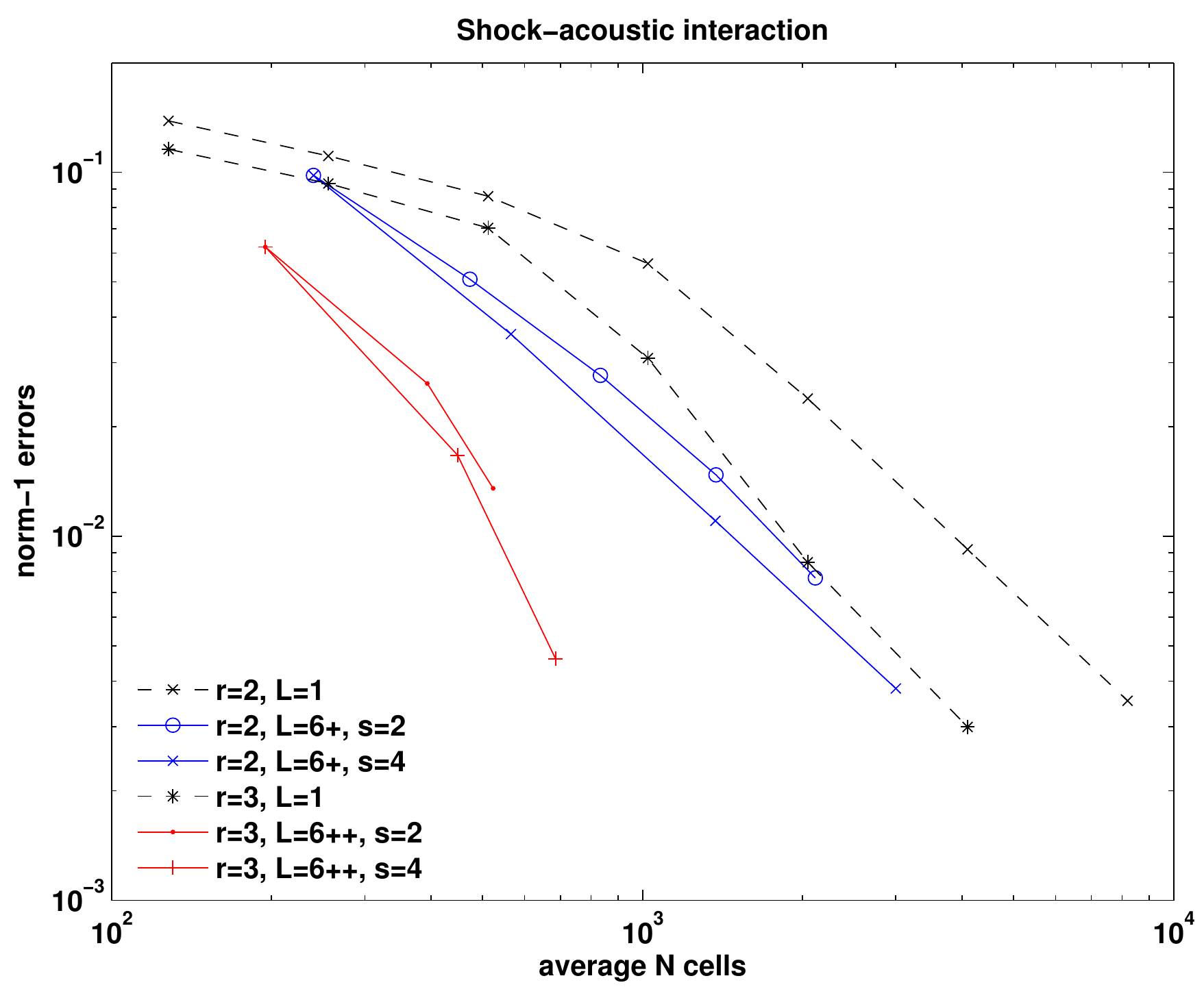}
\hfill
\includegraphics[width=0.49\linewidth]{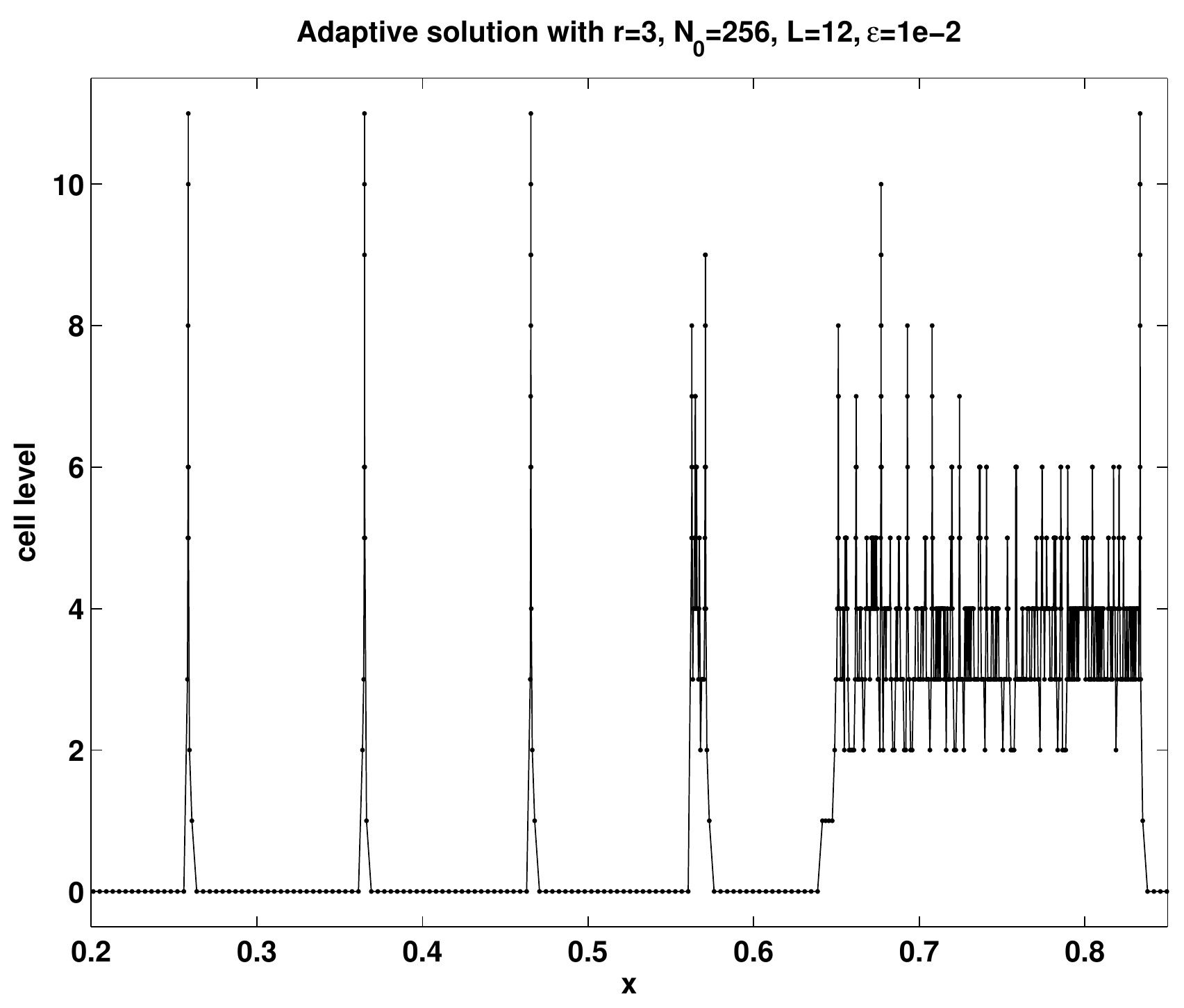}

\includegraphics[width=0.49\linewidth]{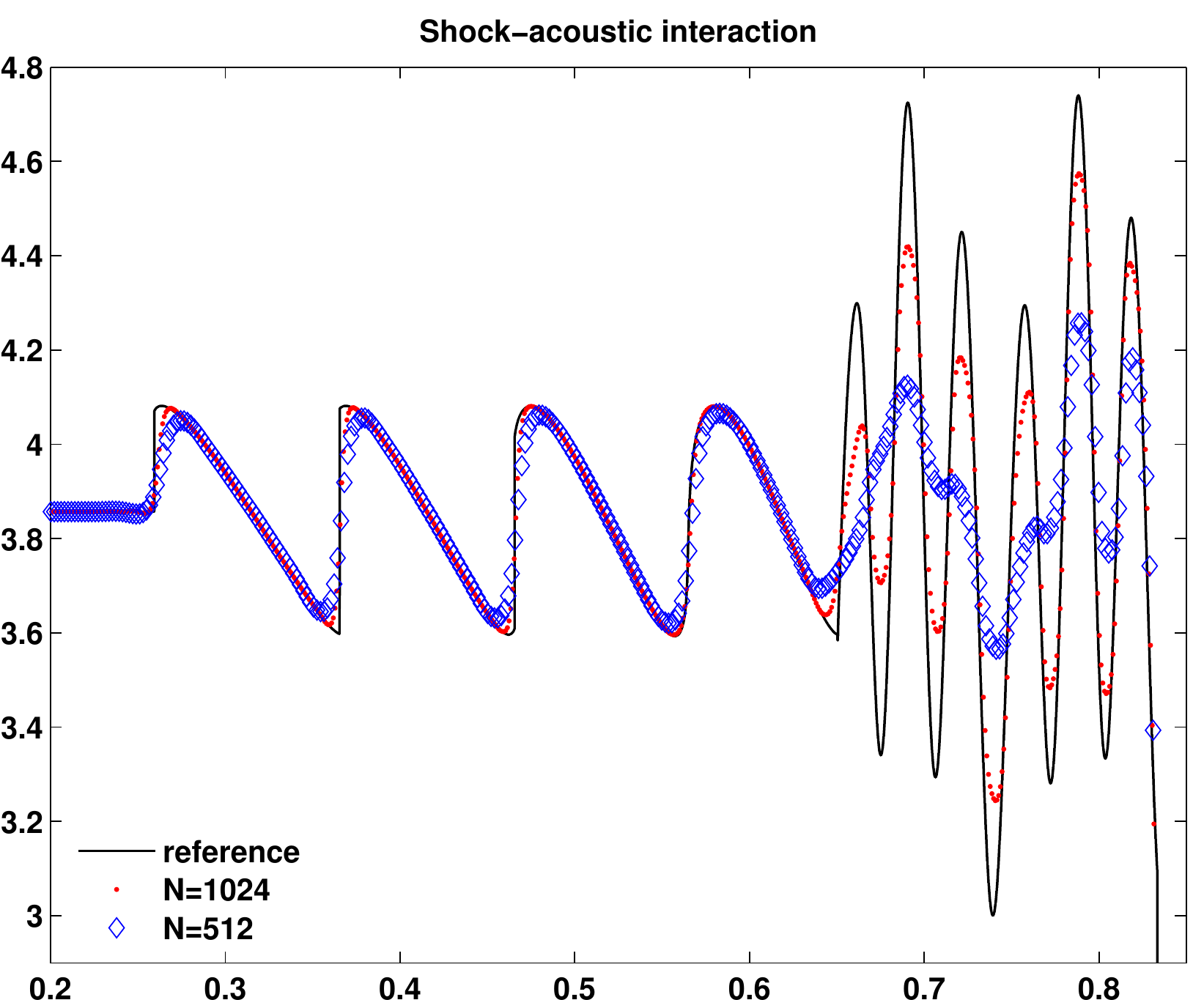}
\hfill
\includegraphics[width=0.49\linewidth]{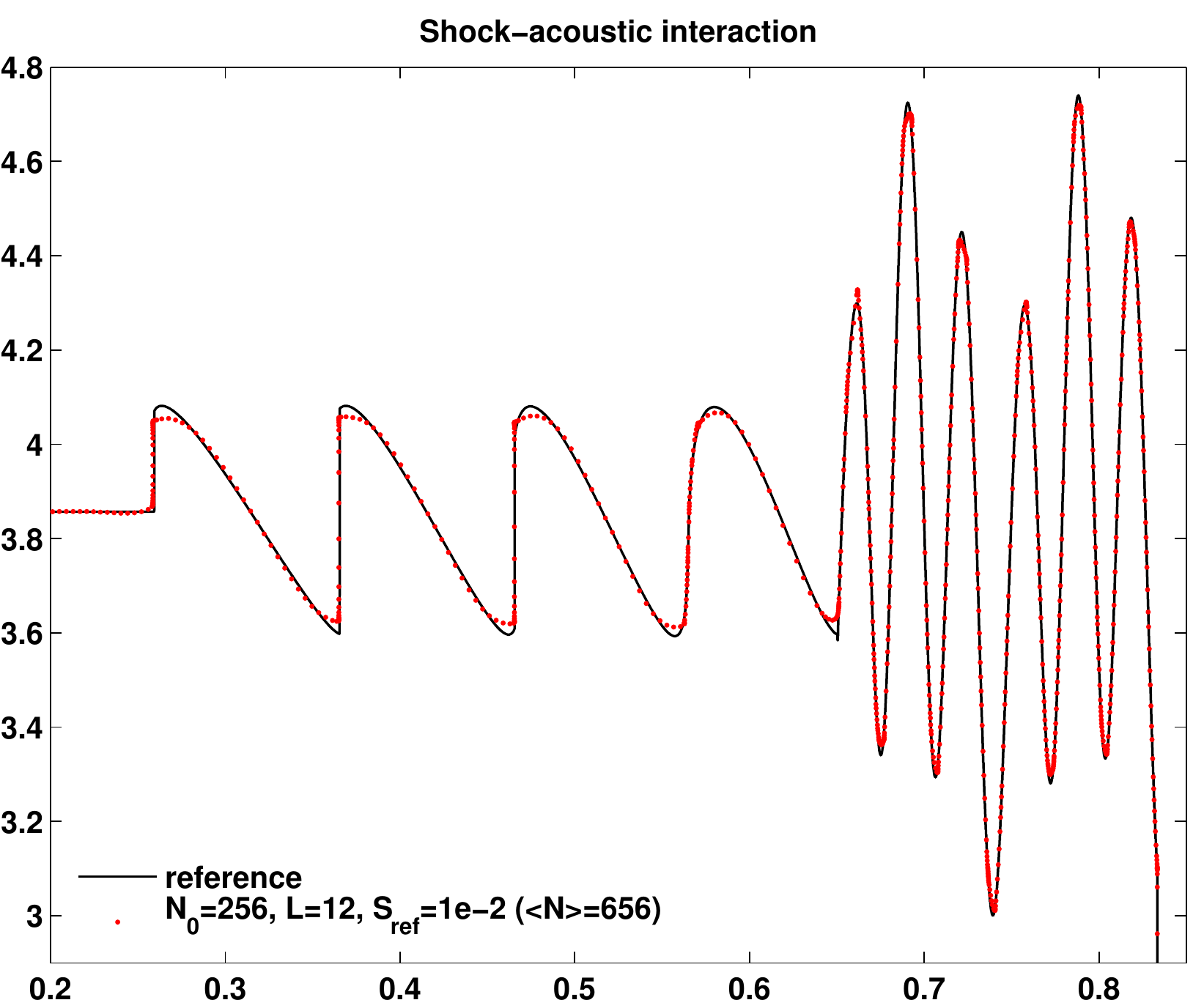}
\end{center}
\caption{Shock-acoustic interaction problem.
(Top-left) Norm-1 errors versus average number of cells, using a reference solution computed by the third order scheme on a uniform mesh of $262144$ cells. The solution on an h-adapted grid is shown with the cell levels (top-right) and with a zoom in the oscilaltory region (Bottom-right). For comparison we show two uniform grid solutions (bottom-left).}
\label{fig:shuOsher}
\end{figure}

\paragraph{One dimensional gas-dynamics: shock-acoustic interaction}
Next we consider the interaction of a shock wave with a standing 
acoustic wave from \cite{SO88}. The conservation law is the one-dimensional Euler
equations with $\gamma=1.4$ and the initial data is
\[ 
(\rho,v,p)=\begin{cases}
  (3.857143,2.629369,10.333333) &, x\in[0,0.25]\\
  (1.0 + 0.2\sin(16\pi x), 0.0, 1.0) &, x\in (0.25,1.0]
\end{cases}
\]
The evolution was computed up to $t=0.2$. As the right moving shock impinges in the stationary wave, a very complicated smooth structure emerges and then gives rise to small shocks and rarefactions (see Fig.\ \ref{fig:shuOsher}). 

The top-left panel of Fig.\ \ref{fig:shuOsher} shows the 1-norm of the errors versus the average number of cells used during the computation.
Due to the high frequency of the waves behind the shock, uniform grid methods starts showing some convergence only when using more than $1000$ cells: before that they lack the spatial resolution required to resolve the smooth waves in the region $x\in[0.65,0.85]$. 

In the adaptive tests, the number of levels is $6$ when the coarse grid has $N_0=32$ cells, so that even the smallest h-adaptive runs have the chance of using cells of size $1/1024$, and we took $S_0=10^{-1}$ for both $r=2,3$.
As in the case of Burgers' equation, in order to observe high convergence rates, one has to increase the number of refinement levels each time that the coarse grid is doubled. Here we additionally observe that the scaling factor $s=4$ gives better results than $s=2$; although, given the presence of the shock, this is somewhat unexpected, it can be understood by looking at the refinement pattern in the top-rigth panel of the Figure, which is representative of the grids employed in the h-adaptive tests. The smallest cells are used only at the shocks, but behind the main shock, there is a non-negligible smooth region where mid-sized cells are employed. These h-adapted grids are an hybrid situation between grids refined only at shocks (which would call for the choice $s=2$) and grids refined in the smooth regions (which would favour the choice $s=8$ in the case of the third order scheme).
Despite the aforementioned difficulties in selecting optimal parameter values, we point out that the second order h-adaptive scheme gives convergence rates similar to those of the uniform grid case, but with lower error constants, while the third order also shows improved error decay rates: the EOC for the two segments shown in the Figure for the run with \texttt{L=6++, s=4} are $1.6$ and $3.0$.

In the lower part of Fig.\ \ref{fig:shuOsher} we compare the reference solution with those computed on uniform grids (left) and with the adaptive algorithm (right), restricted to the interval $x\in[0.2,0.85]$. Solutions computed on uniform grids converge rather slowly, lacking resolution for both the frequency of the waves behind the shock and their amplitude. The adaptive solution captures perfectly
the frequency of the small waves behind the shock and approximates reasonably well their amplitude even with $656$ cells (on average during the time evolution), by making an effective use of adaptivity. In fact the cell sizes are distributed in the computational domain (top-rigth panel) by concentrating them around shocks and high frequency waves, while using larger cells in smooth regions, where the solution is efficiently resolved by the third order scheme.

\paragraph{2D Riemann problem} For a two-dimensional test involving shocks, contact discontinuities and a smooth part, we consider the Euler equation of gas dynamics in two spatial dimensions and set up a Riemann problem in the unit square with initial data given by:
\[
\begin{cases}
\rho=0.8, u=0.1, v=0.1, p=1.0 &x<0.5,y<0.5 \\
\rho=1.0222, u=-0.6179, v=0.1, p=1.0 &x<0.5,y>0.5 \\
\rho=1.0, u=0.1, v=0.8276, p=1.0 &x>0.5,y<0.5 \\
\rho=0.531, u=0.1, v=0.1, p=0.4 &x>0.5,y>0.5
\end{cases}
\]

\begin{figure}
\begin{center}
\includegraphics[width=0.49\textwidth]{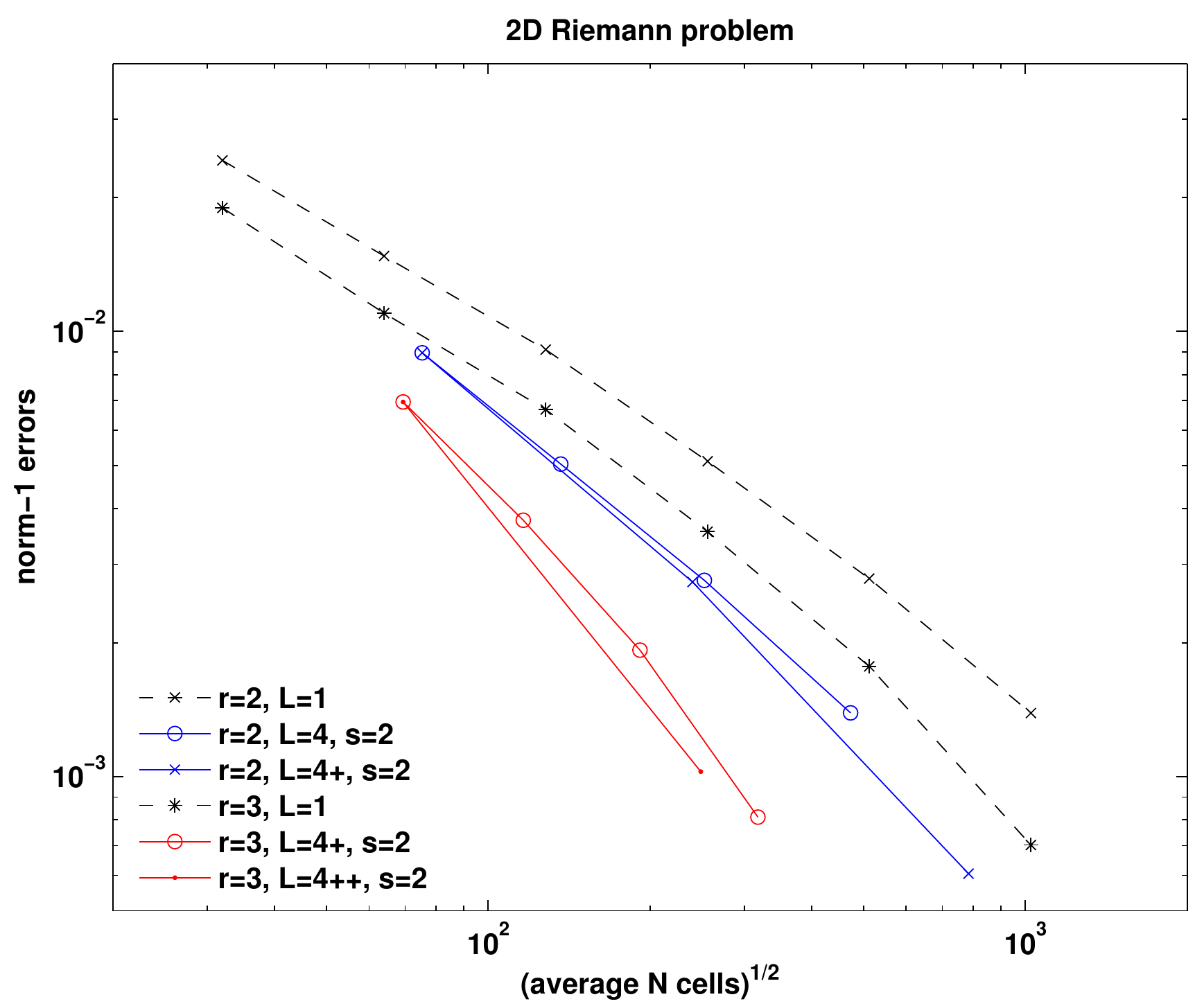}
\includegraphics[width=0.49\textwidth]{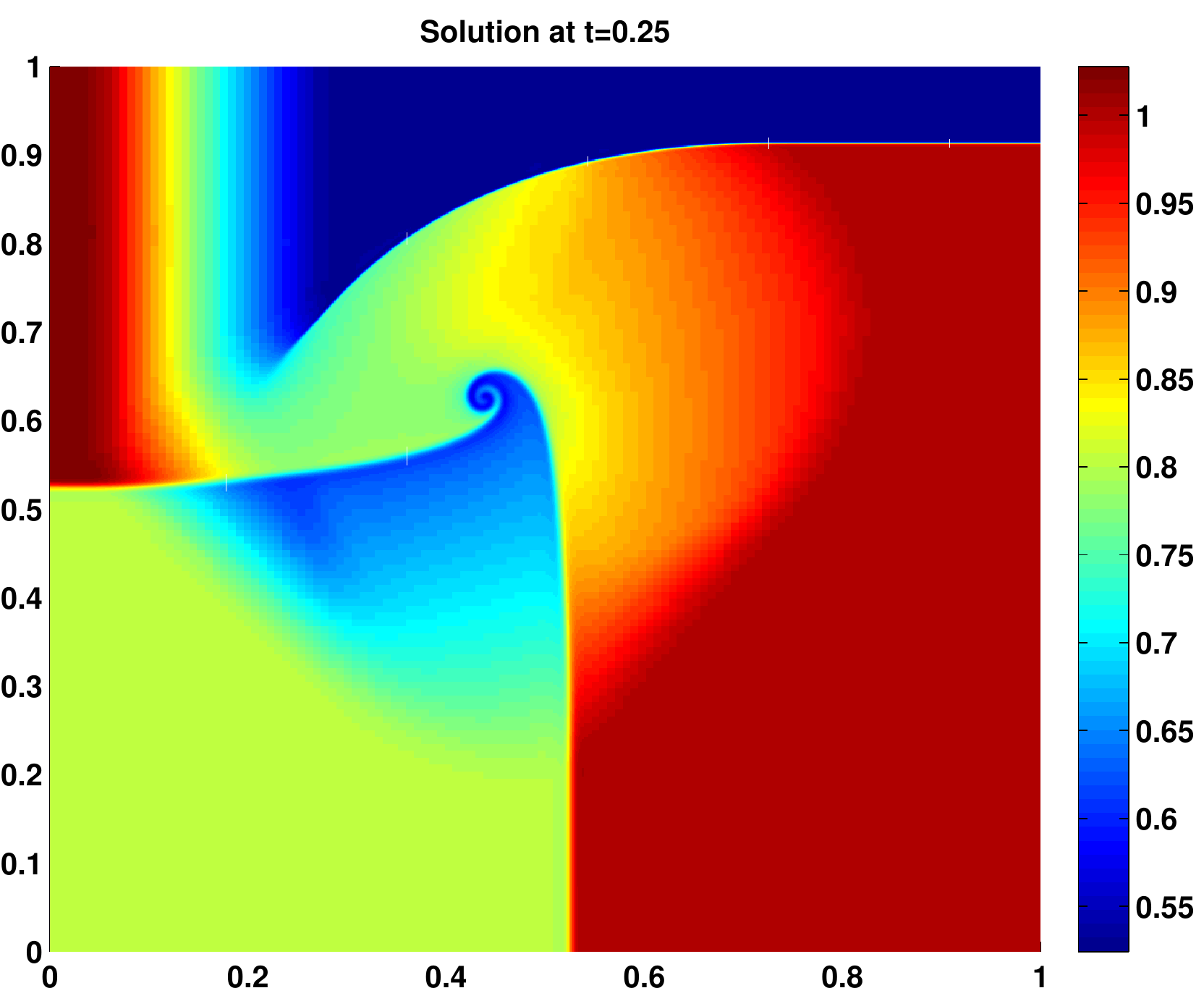}

\includegraphics[width=0.49\textwidth]{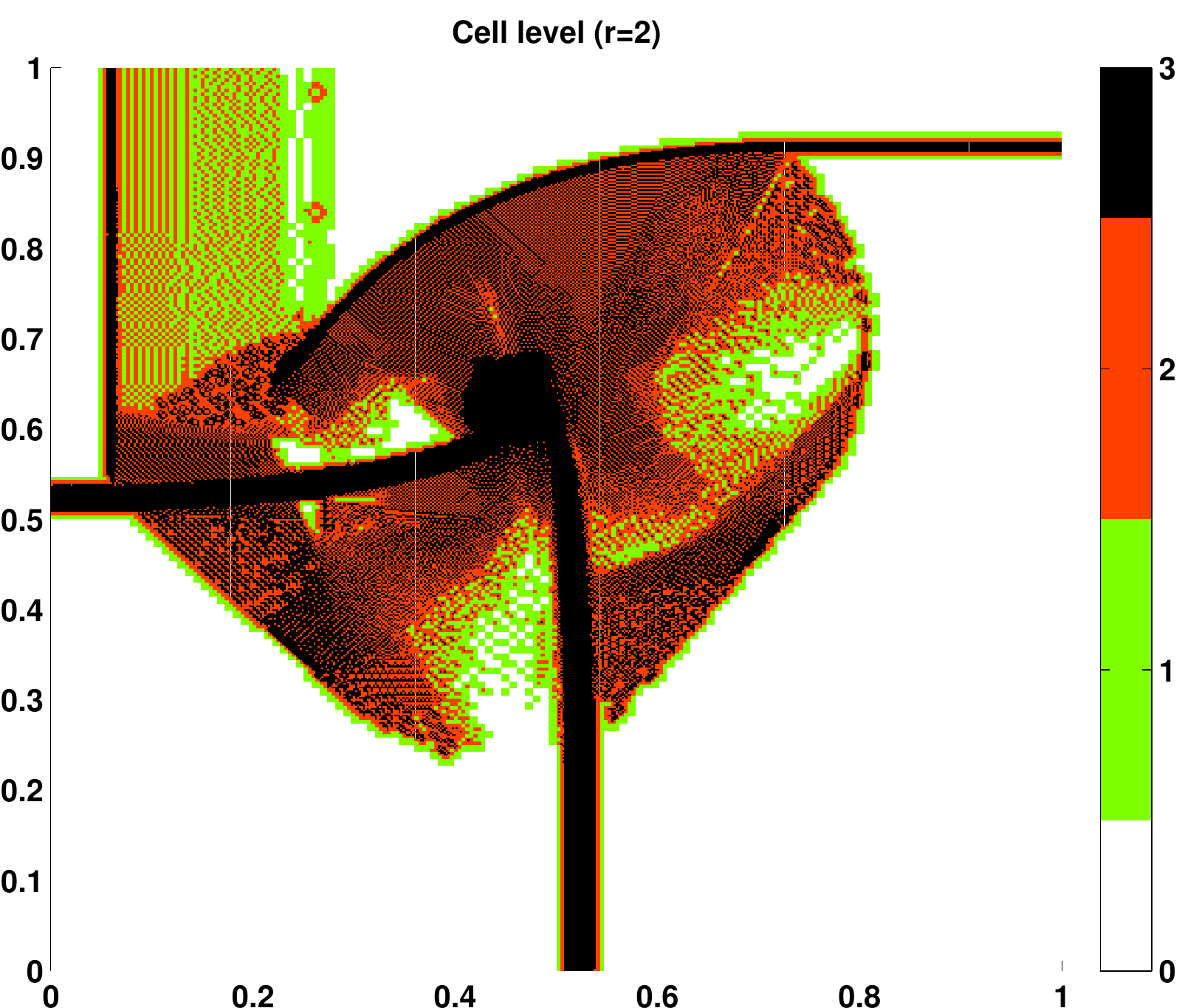}
\includegraphics[width=0.49\textwidth]{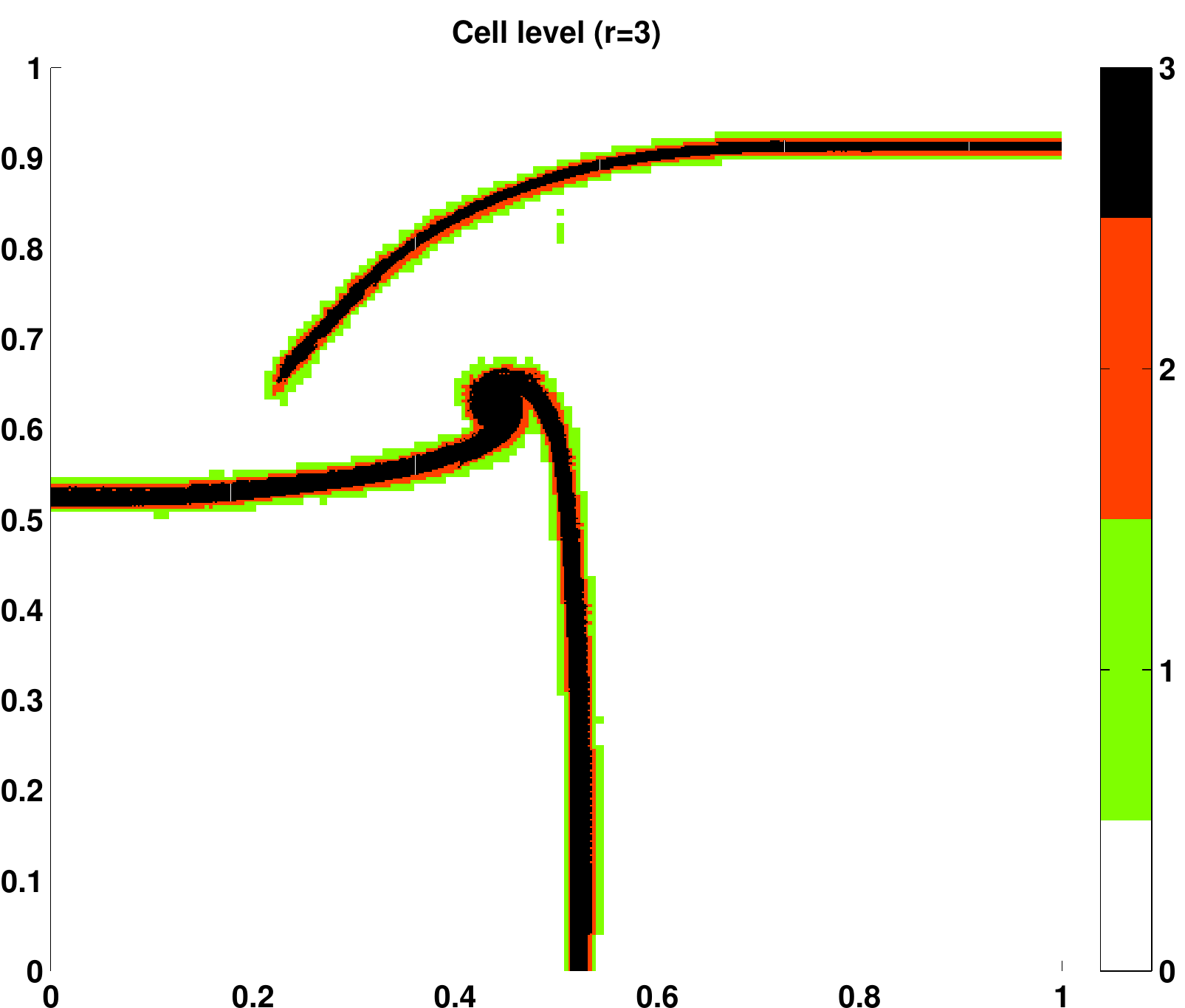}
\end{center}
\caption{Two dimensional Riemann problem for the Euler equations: errors (top-left), density at final time (top-right), grid refinement levels with $N_0=128$, $\mathtt{L}=4$, $S_{\text{ref}}=2\cdot10^{-3}$ (bottom).}
\label{fig:riemann}
\end{figure}

The solutions were compared to a reference one computed with the third order scheme on a uniform grid with $2048\times2048$ cells and the results are presented in Fig.\ \ref{fig:riemann}, whose top-rigth panel depicts the density at final time.
In the error versus number of cells graph (top-left), we can observe that, due to the presence of the shocks, the performance of the two uniform grid schemes have little differences, with the third order one being characterized by approximately the same error decay rate of the second order scheme ($1.0$), only with a better error constant. 
The adaptive schemes ($M=16$, $S_0=10^{-3}$) yield better results that the uniform grid ones, by exploiting small cells to reduce the error in the region close to the shocks, where the numerical scheme are obviously only first order accurate. The EOC are $1.0$ and $1.2$ for the {\tt r=2, L=4+} series and $1.5$ for the {\tt r=3, L=4++} segment. We note however that the third order scheme is much more effective also due to its ability to employ much coarser grids, as it is evident comparing the two computational grids in the lower panels of Fig.\ \ref{fig:riemann}, that refer to the solutions computed with $N_0=128$, $\mathtt{L}=4$ and $S_{\text{ref}}=2\cdot10^{-3}$. As a result the third order scheme can compute a solution with error in the $10^{-3}$ range with $6.18\times10^4$ cells, whereas the uniform grid ones requires $1.05\times10^6$ and the second order adaptive scheme $6.14\times10^5$.

\paragraph{ Shock-bubble interaction in 2D gas-dynamics} For the final test we consider again the 2d Euler equations of gas dynamics and set up an initial datum with a right-moving shock that impinges on a standing bubble of gas at low pressure, as in \cite{CadaTorrilhon09}.
In particular, the domain is $[-0.1,1.6]\times[-0.5,0.5]$ and in the initial datum we distinguish three areas: the left region (A) for $x<0$, the bubble (B) of center $(0.3,0.0)$ and radius $0.2$, the right region (C) of all points with $x>0$ that are not inside the bubble. The initial datum is
\[
\begin{cases}
\rho=3.6666666666666666, u=2.7136021011998722, v=0.0, p=10.0 &(x,y)\in A \\
\rho=0.1, u=v=0.0, p=1.0 &(x,y)\in B \\
\rho=1.0, u=v=0.0, p=1.0 &(x,y)\in C \\
\end{cases}
\]
Boundary conditions are of Dirichlet type on the left, free-flow on the right, reflecting on $y=\pm0.5$.
Due to the symmetry in the $y$ variable, the solution was computed on the half domain with $y>0$, with reflecting boundary conditions on $y=0$.

\begin{figure}
\begin{center}
\includegraphics[width=0.9\textwidth]{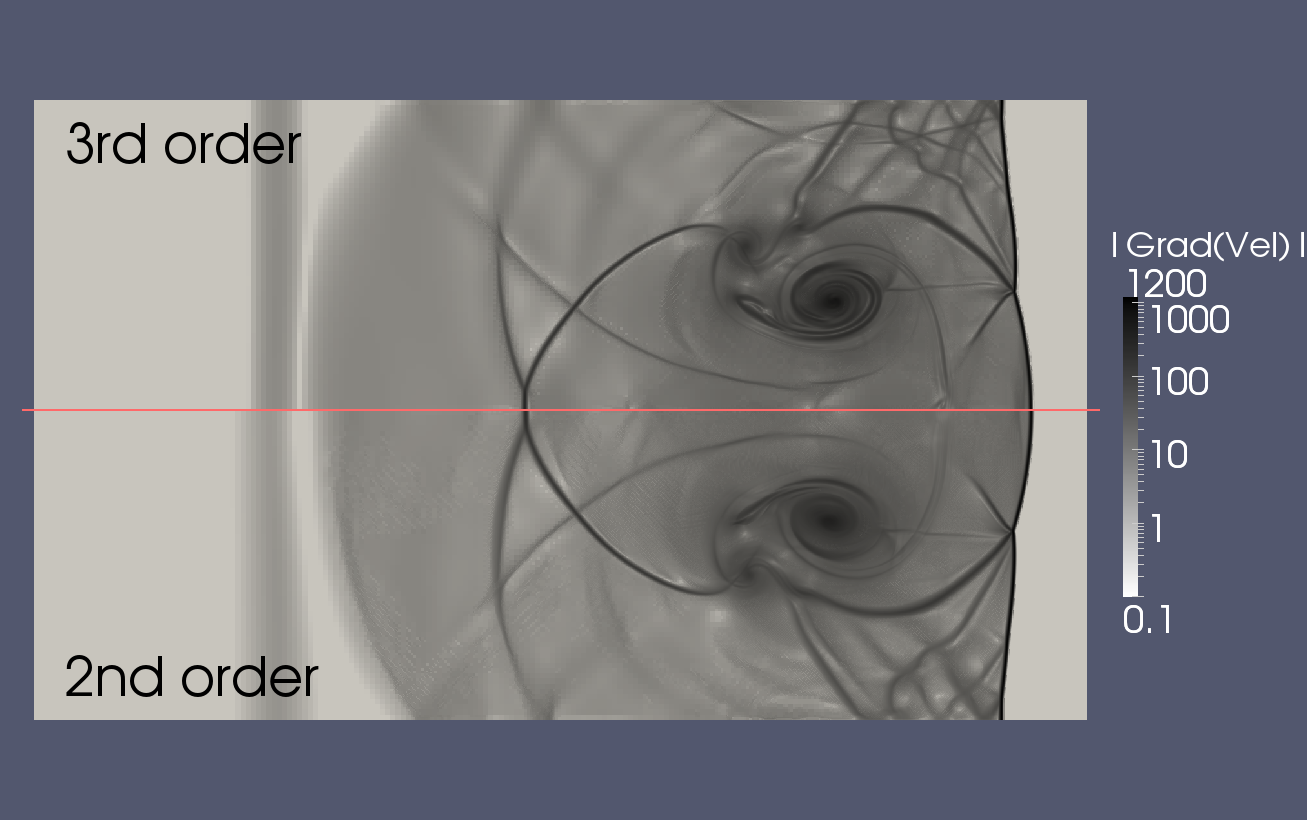}

\includegraphics[width=0.9\textwidth]{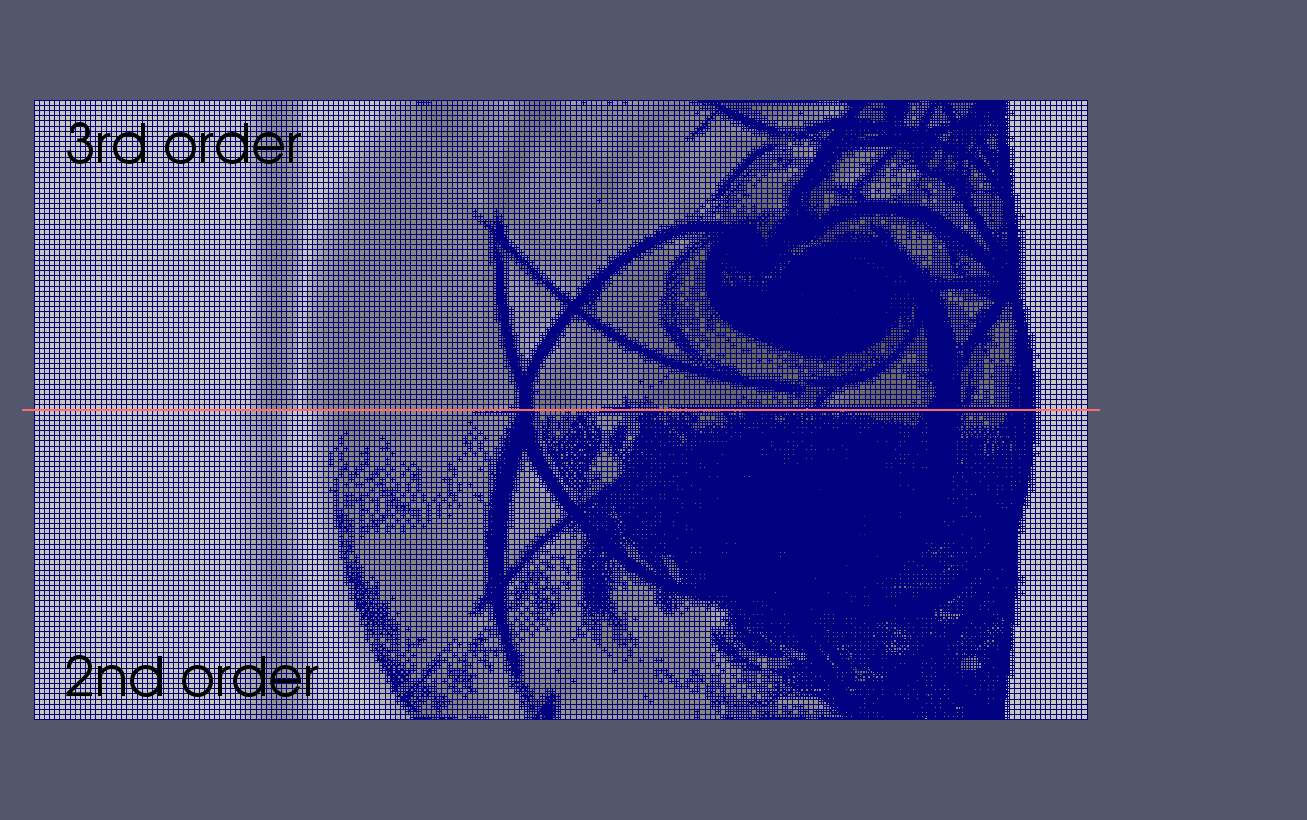}
\end{center}
\caption{Schlieren plots of the solutions (top) of the shock-bubble interaction problem computed with an adaptive grid with $4$ levels from a coarse grid of $204\times60$ cells. The grids at final time are shown in the bottom graph. Since the solution is symmetric, each plot is split in two parts, with the third order solution in the top half and the second order one in the bottom half.}
\label{fig:shockbubble}
\end{figure}

\begin{table}
\begin{tabular}{l|r|r|}
Scheme & norm-1 error & CPU time \\ \hline
$r=2, \mathrm{L}=1, N_0=480\times1632$ & $2.97\cdot10^{-2}$ & $20^h36'$\\
$r=3, \mathrm{L}=1, N_0=240\times816$ & $3.67\cdot10^{-2}$ & $19^h12'$\\
$r=2, \mathrm{L}=4, N_0= 60\times204$ & $2.97\cdot10^{-2}$ & $18^h44'$\\
$r=3, \mathrm{L}=4, N_0= 60\times204$ & $9.60\cdot10^{-3}$ & $18^h35'$\\
\end{tabular}
\caption{Comparison of norm-1 errors and CPU-times for the computation of the shock-bubble interaction problem. These tests were run on a workstation equipped with a 2.60GHz Intel Xeon processor. }
\label{tab:cputimes}
\end{table}

In this tests we used coarse grids with multiples of $15\!\times\!51$ square cells.
Fig.\ \ref{fig:shockbubble} shows a comparison of the solution at $t=0.4$ and the h-adapted computational grid at final time. Each panel depicts a solution in the whole domain, using data from the third order scheme in
the upper half and data from the second order one in the lower half. Both solutions depicted in the Figure have been computed with a coarse grid of $204\times60$ cells, $\mathtt{L}=4$ and $S_{\text{ref}}=2.5\cdot10^{-2}$. The CPU times were comparable ($18^h44^m$ for $r=2$ and $18^h35^m$ for $r=3$), but we can appreciate the much increased sharpness of the solution computed with the third order scheme: note for example the regions close to $y=0.5$ and just behind the main shock, the region of the vortex and the areas just behind it. Also in this case, this better quality solution is obtained using a much coarser grids, i.e. using h-adaptivity in much smaller areas, quite concentrated around the problematic features of the solution. In fact the computation with the second order scheme employed almost twice as much cells than the third order one: $141539$ versus $71895$ cells on average and $2214934$ versus $123114$ at final time. Table \ref{tab:cputimes} show that the computation times compare favourably with the uniform grid schemes. Extra saving in CPU time may be achieved by implementing a local timestepping scheme for time advancement.

\section{Conclusions}
In this work we extended the Compact WENO reconstruction of \cite{LPR:2001} to the case of non-uniform meshes of quad-tree type and used it to build a third order accurate numerical scheme for conservation laws. Indications on how to extend the result to three-dimensional oct-tree grids are provided. 
Even in one space dimension, there are two advantages of the Compact WENO reconstruction over the standard WENO procedure. First, the reconstruction provides a uniformly-accurate reconstruction in the whole cell, which was exploited e.g. in \cite{PS:shentropy} for the integration of source terms, and it is not obtained through dimensional splitting, thus including also the $xy$ terms. Secondly, the linear weights do not depend on the relative size of the neighbors and thus the generalization to two-dimensional grids does not need to distinguish among very many cases of local grid geometry.

The simplicity in computing the reconstruction on non-uniform grids was exploited for coding an h-adaptive scheme that relies on the numerical entropy production as an error indicator. Our tests show that the third order scheme works better than second order ones, not only in the sense that a smaller error is produced in each cell, but also in the sense that, when employed in the h-adaptive scheme, it gives rise to much smaller meshes, that are refined only in a narrow region around shocks.

In order to evaluate the effectiveness of the adaptive schemes we presented some heuristics for the relation between error and number of cells in an h-adaptive scheme in the presence of isolated shocks. In one space dimension we can observe the predicted third order error decay rate, while in the two dimensional case we observe convergence rates slightly higher than those predicted by the heuristics.

Future work of this subject should certainly look at analyzing more deeply the role of the parameter $\epsilon$ appearing in the nonlinear weights of the CWENO reconstruction on non-uniform grids, extending the reconstruction techniques to quad-tree grids composed of triangles split by baricentric subdivision and finally at improving the CPU-time efficiency by a third order local timestepping technique.

\section*{Acknowledgments}
This work was supported by ``National Group for Scientific Computation (GNCS-INDAM)''

%\bibliographystyle{abbrvnat}
%\bibliographystyle{spmpsci}
%\bibliography{crs}

\end{document}